# Predator-Prey Interactions in Communities with prey dispersal and Allee effects

**F. Berezovskaya[1], S. Wirkus[2], C. Castillo-Chavez[3, 4, 5, 6]**

[1]Department of Mathematics, Howard University, Washington D.C., 20059, USA

[2]Division of Mathematical & Natural Sciences, Arizona State University, Glendale, AZ 85306, USA

[3] Mathematics, Computational and Modeling Sciences Center, Arizona State University PO Box 871904, Tempe, AZ, 85287-1904, USA

[4]School of Human Evolution and Social Change, ASU, Tempe, AZ 85287

[5]School of Mathematical and Statistical Sciences, ASU, AZ 85287

[6]Santa Fe Institute, SFI, 1399 Hyde Park Road, Santa Fe, NM, 87501

## Abstract

The population dynamics of predator-prey systems in the presence of patch-specific predators are explored in a setting where the prey population has access to both habitats. The emphasis is in situations where patch-prey abundance drives prey-dispersal between patches, with the fragile prey populations, that is, populations subject to the Alee effect. The resulting four-dimensional model's mathematical analysis is carried out via sub-models that focus in lower dimensional settings. The outcomes depend on, and in fact they are quite sensitive to, the structure of the system, the range of parameter values, and initial conditions. We show that the system can support multi-stability and a diverse set of predator-prey life-history dynamics that includes rather complex dynamical system outcomes. It is argued that in general evolution should favor heterogeneous settings including Allee effects, prey-refuges, and patch-specific predators.

## 1. Introduction.

### 1.1. Background

The *pioneering* work of Lotka and Volterra [23], [29], [30] brought to center stage the importance of developing theoretical frameworks that increase our understanding of the role that predator-prey or competitive or mutualistic interactions have in shaping community structure. This line of theoretical/mathematical research, begun nearly a century ago, continues to challenge and interest ecologists as well conservation and evolutionary biologists. Models incorporating movement within and between sub-populations have been widely investigated in an effort to understand the role of individuals' movement on community sustainability [11], [12], [13], [18], [9], [10], [7].

The study of predator-prey dynamics, broadly understood to include, for example, host-parasite interactions, is of importance in population biology. Theoretical studies that focus on the role of prey-refuges on predator-prey systems have been conducted ([18], [24] and references therein). Post, et al. [26] have focused on the dynamics of two non-interacting prey populations in an environment where the predator switches in response to prey frequency, a response that has a rather strong stabilizing effect on the system. In fact, predators' switching behavior can "control" the system's dynamics, to the point that the predator is able to eliminate the possibility of complex dynamics. Lopez-Gomez et al. [24] have focused on the role of critical patch size on prey survival in systems that do not include predators explicitly. Kuang and Takeuchi [18] have examined the dynamics of predator-prey systems when the prey disperses in response to local (density-dependent) competition showing, for example, that low and high dispersal rates can de-stabilize such systems. Here, we explore the impact of patch specific predators (preference) in a two-patch prey system connected by prey dispersal. The possibility that one of the patches wll serve as a fragile prey refuge (Allee effect, [1]) has been recently analyzed in ([7]). Predator-prey systems where the prey has strong ties to its environment have also been conducted (see [16], [20]).

1.2. **Model Description**

A two-patch model consisting of a predator-prey system with a diffusely migrating prey is the starting point of this manuscript. It is assumed that a fragile prey population (Allee effect, [1], [25]) connects (via its movements) two distinct habitats. We let $u_i \geq 0,\ v_i \geq 0,\ i = 1,2$ denote the population densities of the interacting preys and predators, respectively, in the $i$-th patch. The model's equations are:

$$
\begin{aligned}
u_1' &= \beta_1 f(u_1) - u_1 v_1 + \alpha_1 (u_2 - u_1) \equiv F_1(u_1, v_1, u_2, v_2),\\
v_1' &= \gamma_1 v_1 (u_1 - m_1) \equiv G_1(u_1, v_1, u_2, v_2),\\
u_2' &= \beta_2 f(u_2) - u_2 v_2 + \alpha_2 (u_1 - u_2) \equiv F_2(u_1, v_1, u_2, v_2),\\
v_2' &= \gamma_2 v_2 (u_2 - m_2) \equiv G_2(u_1, v_1, u_2, v_2),\\
&\textit{where} \quad f(u_i) = u_i (u_i - l_i)(1 - u_i)
\end{aligned}
\tag{1}
$$

$\beta_i > 0$ characterize the rates of prey growth; $0 \leq l_i \leq 1$ denote the critical densities of the prey population; $\gamma_i \geq 0$ denote the coefficients of conversion of prey into predator biomass; $m_i \geq 0$ is a measure of the predators' adaptation to the preys; $\alpha_i \geq 0$ characterize migrations of preys in $i$-th patch.

First, we focus on the "symmetric" case:

$$\alpha_1 = \alpha_2 \equiv \alpha,\ \gamma_1 = \gamma_2 \equiv \gamma,\ m_1 = m_2 \equiv m,\ l_1 = l_2 \equiv l,\ \beta_1 = \beta_2 = 1,$$

calling (1s) the symmetric system (1).

1.3. **Overview of Subsystems**

The mathematical analysis of Model (1s) advances through an approach that builds on the analyses of lower dimensional sub-models, see Fig. 1. For example, we first consider the case when prey dispersal is impossible. Thus we have a pair of uncoupled predator-prey (two dimensional) systems. In each patch ( $\alpha = 0$ ) the model describes the dynamics of the densities of prey (u) and predator(v):

$$u' \equiv \frac{du}{dt} = f(u) - uv, \qquad (2)$$
$$v' \equiv \frac{dv}{dt} = \gamma v(u - m)$$

where $f(u) = u(u-l)(1-u)$ and parameters $l$, $m$, $\gamma$ are defined as above.

System (2) modifies and enhances the classical Volterra model. It has been proposed and investigated in prior works (see [5], [8], [28], etc). The phase-parameter portrait of system (2) is shown in Fig.2 and is described below:

**Theorem 1.1.** *For any fixed positive* $\gamma$

1) *and parameters* $(l,m) \in \mathrm{M}\{0 \le l \le 1, m \ge 0\}$ *system (2) in the quadrant* $u \ge 0, v \ge 0$ *has equilibria* $O(0,0), O_1(1,0), O_l(l,0)$ *and equilibrium* $A(m,(m-l)(1-m))$ *if* $0 \le l < m$.

2) *the parameter space* $\mathrm{M}$ *is divided into 5 regions of qualitatively different phase portraits of System (2). Boundaries between regions correspond to the bifurcations of co-dimension* 1:

***S*₁**: $m = 1$ *and* ***S****l*: $m = l$, *the appearance/disappearance of point* $A$ *in the* 1st *quadrant by transcritical bifurcations with* $O_1$, $O_l$, *respectively;*

$\boldsymbol{H_1}: m = \dfrac{l+1}{2}$, *the change of stability of* $A$ *in the supercritical Andronov-Hopf bifurcation (with appearance/disappearance of a stable limit cycle);*

***L:*** *m=m (l) the disappearance/ appearance of a stable limit cycle in a heteroclinics composed by the separatrices of equilibria* $O_1$ *and* $O_L$.[1]

Thus, the local (one-patch) Model (2) demonstrates the possibility of prey-predator coexistence in stable equilibrium or in oscillations if parameters are in the parameter domains 2 and 3 of Fig. 2 (see also 3.1). We show in this work that prey dispersal between two bilocal (two-patch) models (1) "originates" new dynamical modes of population coexistence and essentially generalizes the possibility of population persistence.

The previous case when prey dispersal is impossible –System (2)—sets the stage for the study of the predator-prey dynamics via the invasion of patch-specific predators. Hence, we also look at the impact of dispersal on the dynamics of a two-patch predator-free environment. Let us note that system (1s), which was developed as a "two-dimension $\alpha-$updating" of Model (2) can be considered also as a "two-dimension updating" of the model (see, [12], [13]):

$$\begin{aligned} u_1' &= f(u_1) + \alpha(u_2 - u_1) = F_1(u_1, 0, u_2, 0), \\ u_2' &= f(u_2) + \alpha(u_1 - u_2) = F_2(u_1, 0, u_2, 0), \end{aligned} \qquad (3)$$

which describes the dynamics of two Allee-type prey populations interacting diffusely. From this point of view System (1s) determines the role of predators in a population community. More exactly, it models how the behavior of the community Model (3) should change with the introduction of their predators. As we will show later, dispersal in a predator-free two-patch environment where the prey must reach a critical mass to survive can support up to nine equilibria including five boundary (one or two prey populations are absent) and four co-internal positive equilibria (both prey populations are present) when the rate of dispersal is low.

[1] Curve *m(l)* was found numerically in [5]; it was recalculated with help of the specific computer algorithm in [28]

The next subsystem that we consider is the situation when one patch faces predation while the other is a refuge (no access to predators). Model (1) can be thought up as a "one-dimension updating" of a two-patch models [7].

$$
\begin{aligned}
u_1' &= f(u_1) + \alpha(u_2 - u_1) = F_1(u_1, 0, u_2, v_2), \\
u_2' &= f(u_2) - u_2 v_2 + \alpha(u_1 - u_2) = F_2(u_1, 0, u_2, v_2), \qquad \text{(4a)} \\
v_2' &= \gamma v_2 (u_2 - m) = G_2(u_1, 0, u_2, v_2)
\end{aligned}
$$

and

$$
\begin{aligned}
u_1' &= f(u_1) - u_1 v_1 + \alpha(u_2 - u_1) = F_1(u_1, v_1, u_2, 0), \\
v_1' &= \gamma v_1 (u_1 - m) = G_1(u_1, v_1, u_2, 0) \qquad \text{(4b)} \\
u_2' &= f(u_2) + \alpha(u_1 - u_2) = F_2(u_1, v_1, u_2, 0),
\end{aligned}
$$

which describe the dynamics of communities consisting of prey and predator when the prey can disperse between both patches. Systems (4a) and (4b) differ only by designation of variables; we thus omit indices and refer to it as System (4). As we will show in later sections, when predators have access to one patch (there is a predator-free or a prey-refuge) the system will support one to three positive equilibria (both prey populations and the predator surviving). These three-dimensional positive equilibrium points correspond to boundary (pre-prey plane) equilibria in the absence of the predator.

We show below that the dynamics of System (1s) includes the dynamics of (2), (3) and (4). In other words, the two-dimensional prey-prey system is naturally embedded in the refuge-prey-predator system. Similarly, the refuge-prey-predator system is also naturally embedded in the full predator-prey-prey-predator system. Additionally, however, system (1s) produces its "own" non-trivial 4D behaviors, including the equilibrium $AA(m,(m-l)(1-m),m,(m-l)(1-m))$ and additional oscillations. These behaviors are of the *main attention* in this work. We will distinguish between "*trivial*" and "*non-trivial*" equilibria, the former having at least one zero-coordinate. A trivial equilibrium of Systems (1s) and (4) can arise from non-trivial equilibria of system (2) or system (3) and inherit some properties of the lower dimension equilibria. In an effort to understand the complete dynamics of the model community close to the non-trivial point *AA* we will analyze all equilibria of the model and other non-trivial modes.

Some of the equilibria have their coordinates given exactly while others have an asymptotic expansion (explicitly to $O(\alpha)$ or $O(\alpha^2)$, see Table 1). We show that for $0 < l < m \le 1$ the system can have up to *sixteen* equilibria from which we distinguish "*strictly symmetric*" equilibria: $u_1 = u_2$, $v_1 = v_2$, "*three-dimensional*" equilibria possessing one zero-coordinate: $v_1 = 0$ or $v_2 = 0$, and "*two-dimensional*" equilibria possessing two zero-coordinates: $v_1 = v_2 = 0$.

The paper is organized as follows. In Sections 2 and 3 we analyze positions and stability of equilibrium points of model Systems (1s) and its subsystems. We consider also the problem of the change of stability with increasing dimension of the subsystems. The results of this analysis in asymptotic (on $\alpha$) form are done in Tables 1,2. In Section 4 we describe the global dynamics of the subsystems in phase-parameter spaces. Section 5 contains the results of the analytical and computer analysis of 4D-behaviors of the model (1s) that are presented in the form of phase-parameter portraits. Discussion of the outcomes and their biological interpretations are done in Section 6.

## 2. Coordinates of equilibria of models (1s), (2), (3), (4)

### 2.1. *Coordinates of 2D-equilibria*

For $\alpha = 0$ as well as for $u_1 = u_2$ the Model (1s) describes two independent subsystems in form (2). System (2) has *trivial* equilibria $O_0(0,0), O_l(l,0), O_1(1,0)$ and *non-trivial* equilibrium $A(m, v^* = (m-l)(1-m))$ in the 1st quadrant if $0 < l < m < 1$ (see Fig.1).

For any $\alpha$ System (3) has trivial equilibrium $O_{00}(0,0)$ and symmetric non-trivial equilibria $O_{ll}(l,l)$, $O_{11}(1,1)$. It can also have up to three pairs of *non-trivial* equilibria $C_1(u_1^*, u_2^*)$ and $C_2(u_2^*, u_1^*)$ where $u_1^*, u_2^*$ are different than the $0, l, 1$ roots of the system

$$\begin{aligned} F_1(u_1,0,u_2,0) &\equiv f(u_1) + \alpha(u_2 - u_1) = 0, \\ F_2(u_1,0,u_2,0) &\equiv f(u_2) + \alpha(u_1 - u_2) = 0. \end{aligned} \tag{5}$$

To specify mutual placing of the equilibria we find equations of the coalescing of equilibrium points $C$ in phase-parameter space. The condition is defined by System (5) with the additional requirement that

$$\left| \frac{\partial(F_1(u_1,0,u_2,0), F_2(u_1,0,u_2,0))}{\partial(u_1,u_2)} \right| \equiv f_u(u_1) f_u(u_2) - \alpha(f_u(u_1) + f_u(u_2)) = 0 \tag{6}$$

In parameter space of Model (3) System (5), (6) defines, by the implicit form, the boundary ***SC,*** which divides domains where the model has no non-symmetric non-trivial equilibria, has one pair and three pairs of those (see Fig. 3). Analysis of the equilibria as well as the boundaries (5), (6) was done by expanding the functions $F_1(u_1,0,u_2,0), F_2(u_1,0,u_2,0)$ in series in $\alpha$ and considering the asymptotics explicitly to $O(\alpha)$. Note, that equilibria $C(u_1^*, u_2^*)$ can be considered as those arising from two of the equilibria $O_0, O_l, O_1$ of System (2) affected by the parameter $\alpha$. Due to this fact we denote the *C*-equilibria as: $C_{0l}, C_{l0}, C_{01}, C_{10}, C_{l1}, C_{1l}$ supposing that $C_{00} \equiv O_{00}$, $C_{ll} \equiv O_{ll}$, $C_{11} \equiv O_{11}$.

The result of the analysis is collected in

***Proposition 2.1.*** *For positive values of parameters the boundary **SC** consists of two branches whose asymptotic (on $\alpha$) forms are:*

$$SC_1 : \alpha_1(l) = 1 - l + l^2 - \frac{[(2l-1)(l-2)(l+1)]^{2/3}}{2},$$

$$SC_2 : \alpha_2(l) = \frac{l - l^2}{2}.$$

*Three pair of equilibria* $C_{01}, C_{10}$, $C_{0l}, C_{l0}$, $C_{l1}, C_{1l}$ *are in the domain bounded by curves:* $\alpha = 0$ *and* $\alpha = \alpha_1(l)$*; one pair of equilibria* $C_{0l}, C_{l0}$ *are in the domain bounded by curves:* $\alpha = \alpha_1(l)$ *and* $\alpha = \alpha_2(l)$*; no equilibria if* $\alpha > \alpha_2(l)$.

*The asymptotic coordinates of the non-symmetric* $C$*-equilibria (explicitly to* $O(\alpha)$*) are given in Table 1.*

2.2. ***Coordinates of 3D- and 4D- equilibria***

Systems (4a) and (4b)

(1) have trivial equilibria $C_1(u_1^*, u_2^*, 0)$ and $C_2(u_2^*, 0, u_1^*)$ where $u_1^*, u_2^*$ are roots of the system (5). The description of these equilibria is done in Proposition 2.1;

(2) have from one up to three pairs of *non-trivial* equilibria $B_1(m, z, y)$ and $B_2(y, m, z)$, correspondingly, where $y, z$ satisfy

$$\begin{aligned} &f(y) + \alpha(m - y) = 0, \\ &z = (f(m) + \alpha(y - m))/m. \end{aligned} \qquad (7)$$

The condition of the coalescing of pairs of equilibrium points $B$ is defined by system (7) with the additional requirement that

$$f_y(y) - \alpha = 0. \qquad (8)$$

In parameter space of Model (4) System (7), (8) defines, by the implicit form, the boundary ***SB*** divided domains where the model has one and three pairs of three-dimensional equilibria (see Fig. 4).

Note, that equilibria $B_1(y,m,z)$, $B_2(m,z,y)$ can be considered as those "arising" from two equilibria of System (2): $O_0, A,\ O_l, A,$ or $O_1, A$. Thus, when systems have three equilibria we call them $B_1^0,\ B_1^l,\ B_1^1$ and $B_2^0,\ B_2^l,\ B_2^1$, respectively, whereas we call them $B_1, B_2$ when systems have only one equilibrium. As above, the analysis of coordinates of the equilibria as well as boundaries (7), (8) was done by expanding the functions $F_1(u_1,0,u_2,v_2), F_2(u_1,v_1,u_2,0)$ in series in $\alpha$ and considering the asymptotics explicitly to $O(\alpha)$ or $O(\alpha^2)$.

The result of the analysis is collected in

***Proposition 2.2.*** *The parameter boundary is the surface* ***SB*** *given by the equation:*

$[27\alpha m - 9(\alpha+1)(1+l) + 2(1+l)^3]^2 + 4(3\alpha - 1 + l - l^2)^3 = 0$ *where* $\alpha < (1-l+l^2)/3$.

*It consists of two branches* $SB_{12}, SB_{23}$

$$SB_{12} : m(\alpha,l) = \frac{9(\alpha+l)(1+l) + 2(1-l+l^2-3\alpha)^{3/2}}{27\alpha},$$
$$SB_{23} : m(\alpha,l) = \frac{9(\alpha+l)(1+l) - 2(1-l+l^2-3\alpha)^{3/2}}{27\alpha}$$

*Corresponding to the fold bifurcation in the system (4) everywhere except on the line* $\alpha = (1-l+l^2)/3,\ m = (1+l)^3/(27\alpha)$ *which corresponds to the cusp bifurcation. There are two equilibria of systems (4a) and (4b) on* $SB_{12} \cup SB_{23}$ *(one equilibrium is one-multiple, the other has multiplicity 2) and the triple equilibria* $B_1(m,z,y)$ *and* $B_2(y,m,z)$ *with* $z = (1+l)/3,\ \ y = (f(m) + \alpha(z-m))/m$ *at the cusp line.*

*System (4) has three pair of equilibria* $B_1^0$, $B_1^l$, $B_1^1$ and $B_2^0$, $B_2^l$, $B_2^1$ *inside the parameter domain bounded by* ***SB*** *and only one pair* $B_1$, $B_2$ *outside this domain (see Fig.4). The asymptotic coordinates of the* $B$*-equilibria exactly to* $O(\alpha)$ or $O(\alpha^2)$ *are given in Table 2.*

Now we are ready to describe coordinates of equilibria of the 4D system (1s).

**Theorem 2.1.**

(1) For arbitrary positive $\alpha$ System (1s)

●has up to three pairs of *trivial "two-dimensional"* equilibria $C_1(u_1*,0,u_2*,0)$ and $C_2(u_2*,0,u_1*,0)$ where $u_1*,u_2*$ *are distinguished from* $0,l,1$ *roots of the system (5); mutual placing of these equilibria is done in Proposition 2.1 (see also Fig.3);*

●*has from one up to three pairs of trivial "three-dimensional" equilibria* $B_1(m,z,y,0,)$ *and* $B_2(y,0,m,z)$ *where* $y,z$ *satisfy (7); mutual placing of these equilibria is done in Proposition 2.2 (see also Fig.4);*

●*has trivial symmetric equilibria* $O(0,0,0,0)$, $O^{ll}(l,0,l,0)$, $O^{11}(1,0,1,0)$ *and if* $0 < l < m < 1$ *has also non-trivial one* $AA(m,v*,m,v*)$ *where* $v* = f(m)/m = (m-l)(1-m)$.

## 3. ***Linear stability of equilibria***

### 3.1. ***On lower dimension non-trivial equilibrium and higher dimension trivial equilibrium***

Equilibria of Model (2) were completely studied in [5], [4], model (3) in [12] and model (4) in [7]. The main attention below will be given to the analysis of the stability of their "images" into the frame of Model (1s). Thus, we will consider points $C_1(u_1*,0,u_2*,0)$ and $C_2(u_2*,0,u_1*,0)$ where $u_1*,u_2*$ are distinguished from $0,l,1$ roots of the system (5), $B_1(m,z,y,0,)$ and $B_2(y,0,m,z)$ where $y,z$ satisfy (7), as well as points $O(0,0,0,0)$, $O^{ll}(l,0,l,0)$, $O^{11}(1,0,1,0)$ and $AA(m,v^*,m,v^*)$ where $v^* = f(m)/m = (m-l)(1-m)$.

(1) We first consider the equilibria denoted as $C$ which are non-trivial in the 2D system (3) and become trivial in both the 3D system (4) and the 4D system (1s).

Jacobians of system (3) at equilibrium $C(u_1^*,u_2^*)$,system (4b) at $C(u_1^*,u_2^*,0)$ and system (1s) at $C(u_1^*,0,u_2^*,0)$ are correspondingly:

$$J_2=\begin{pmatrix} P_1(C) & \alpha \\ \alpha & P_2(C) \end{pmatrix},\ J_3=\begin{pmatrix} P_1(C) & \alpha & 0 \\ \alpha & P_2(C) & 0 \\ 0 & 0 & S_2(C) \end{pmatrix},\ J_4=\begin{pmatrix} P_1(C) & -u_1 & \alpha & 0 \\ 0 & S_1(C) & 0 & 0 \\ \alpha & 0 & P_2(C) & -u_2 \\ 0 & 0 & 0 & S_2(C) \end{pmatrix},$$

where $P_1(C) = f_u(u_1^*) - \alpha,\ P_2(C) = f_u(u_2^*) - \alpha;\ S_1(C) = -\gamma(m-u_1^*),\ S_2(C) = -\gamma(m-u_2^*)$ and $f(u) = u(u-l)(1-u)$.

These Jacobians have characteristic polynomials, which are equal respectively to

$$\phi_2(\lambda) = (P_1(C)-\lambda)(P_2(C)-\lambda) - \alpha^2,\ \phi_3(\lambda) \equiv \phi_2(\lambda)(S_2(C)-\lambda),\ \phi_4(\lambda) \equiv \phi_3(\lambda)(S_1(C)-\lambda). \quad (9)$$

So, each successive characteristic polynomial differs from the previous one only by one factor.

Similar arguments hold for the three-dimensional and four-dimensional points $B$. The Jacobians of systems (4) at $B(y,m,z)$ and (1s) at $B(y,0,m,z)$, respectively,

$$\tilde{J}_3 = \begin{pmatrix} P_1(B) & \alpha & 0 \\ \alpha & P_2(B) & -m \\ 0 & \gamma z & 0 \end{pmatrix},\quad \tilde{J}_4 = \begin{pmatrix} P_1(B) & -y & \alpha & 0 \\ 0 & S_1(B) & 0 & 0 \\ \alpha & 0 & P_2(B) & -m \\ 0 & 0 & \gamma z & 0 \end{pmatrix}$$

where $P_1(B) = f_u(y) - \alpha,\ P_2(C) = f_u(m) - z - \alpha;\ S_1(C) = -\gamma(m-y)$.

These Jacobians have characteristic polynomials, which are equal to

$\varphi_3(\lambda) = -\lambda(P_1(B)-\lambda)(P_2(B)-\lambda) + \gamma mz(P_1(B)-\lambda) + \alpha^2\lambda,\ \ \varphi_4(\lambda) = \varphi_3(\lambda)(S_1(B)-\lambda)$, (10) and so distinguished only by one factor. Thus, the following statement holds

***Proposition 3.1.*** (1) *Let* $u_1, u_2$ *be roots of system (5). Then systems (3), (4) and (1s) have two identical eigenvalues* $\lambda_1,\ \lambda_2$ *around equilibria* $C(u_1,u_2)$, $C(u_1,u_2,0)$ *and* $C(u_1,0,u_2,0)$: $\lambda_i(C(u_1,u_2)) = \lambda_i(C(u_1,u_2,0) = \lambda_i(C(u_1,0,u_2,0),\ \ i=1,2$; *systems (4) and (1s) have three identical eigenvalues* $\lambda_1,\ \lambda_2,\ \lambda_3$ *around equilibria* $C(u_1,u_2,0)$ *and* $C(u_1,0,u_2,0)$: $\lambda_i(C(u_1,u_2,0) = \lambda_i(C(u_1,0,u_2,0),\ \ i=1,2,3$;

(2) *Let* $u_1, v_2$ *be roots of System (7). Then Systems (4b) and (1s) have three identical eigenvalues* $\mu_1,\ \mu_2,\ \mu_3$ *around equilibria* $B(u_1,m,v_2)$ *and* $B(u_1,0,m,v_2)$, *correspondingly*: $\mu_i(B(u_1,m,v_2) = \mu_i(B(u_1,0,m,v_2),\ \ i=1,2,3$. *The same is true for eigenvalues of equilibria B of Systems (4a) and (1s).*

3.2. ***Estimation of eigenvalues***

The arbitrary equilibrium Jacobian matrix $\dfrac{\partial(F_1,G_1,F_2,G_2)}{\partial(u_1,v_1,u_2,v_2)}$ of System (1) has the specific block-diagonal form:

$$J(u_1,v_1,u_2,v_2) = \begin{pmatrix} P_1 & Q_1 & \alpha & 0 \\ R_1 & S_1 & 0 & 0 \\ \alpha & 0 & P_2 & Q_2 \\ 0 & 0 & R_2 & S_2 \end{pmatrix}$$

where

$P_1 = f_u(u_1) - v_1 - \alpha,\ Q_1 = -u_1,\ R_1 = \gamma v_1,\ S_1 = -\gamma(m - u_1),$
$P_2 = f_u(u_2) - v_2 - \alpha,\ Q_2 = -u_2,\ R_2 = \gamma v_2,\ S_2 = -\gamma(m - u_2).$

It can be verified that for certain important cases the characteristic polynomial $\Psi(\lambda)$ associated with $J$ has eigenvalues that are given by the following statement

***Lemma* 3.1.**

*1) If $J$ is block- symmetric matrix, i.e., $P_1 = P_2$, $Q_1 = Q_2$, $R_1 = R_2, S_1 = S_2$, then characteristic polynomial $\Psi(\lambda)$ has real roots:*

$$\lambda_{1,2} = \frac{P_1 + S_1 - \alpha + \sqrt{\Delta}}{2},\ \lambda_{3,4} = \frac{P_1 + S_1 - \alpha - \sqrt{\Delta}}{2}, \qquad (11)$$
$$\Delta = (P_1 - S_1 + \alpha)^2 + 4Q_1R_1;$$

*2) If $R_1 = R_2 = 0$ then characteristic polynomial $\Psi(\lambda)$ has four real roots:*

$$\lambda_1 = S_1,\ \lambda_2 = S_2,\ \lambda_{3,4} = \frac{P_1 + P_2 \pm \sqrt{(P_1 - P_2)^2 + 4\alpha^2}}{2}; \qquad (12)$$

*3) If $Q_2R_2 = 0$ then $\Psi(\lambda)$ has root $\lambda_1 = S_2$ and the other roots satisfy the equation*

$$\lambda^3 - (P_1 + P_2 + S_1)\lambda^2 + (P_1P_2 + (P_1 + P_2)S_1 - Q_1R_1 - \alpha^2)\lambda - (P_1P_2S_1 - P_1P_2S_2 - \alpha^2 S_1{}^2) = 0 \qquad (13)$$

4) *If $Q_1R_1 = 0$ then $\Psi(\lambda)$ has root $\lambda_1 = S_1$ whereas the other roots satisfy the equation*

$$\lambda^3 - (P_1 + P_2 + S_2)\lambda^2 + (P_1P_2 + (P_1 + P_2)S_2 - Q_2R_2 - \alpha^2)\lambda - (P_1P_2S_2 - P_1P_2S_2 - \alpha^2 S_2{}^2) = 0 \qquad (14)$$

Applying Statements 1 and 2 of Lemma 2.1 to system (1s) we get the following description of eigenvalues around symmetric equilibria.

***Proposition 3.2.*** *The eigenvalues of the System (1s) at the symmetric equilibria* $O(0,0,0,0)$, $O^{ll}(l,0,l,0)$, $O^{11}(1,0,1,0)$, $AA(m,v^*,m,v^*)$ *where* $v^* = f(m)/m$ and $0 < l < m < 1$ *are, respectively,*

$$\lambda_1(O) = \lambda_2(O) = -\gamma m,\ \lambda_3(O) = -l,\ \lambda_4(O\ ) = -l - 2\alpha;$$
$$\lambda_1(O^{ll}) = \lambda_2(O^{ll}) = -\gamma(m-l),\ \lambda_3(O^{ll}) = l(1-l),\ \lambda_4(O^{ll}) = l(1-l) - 2\alpha;$$
$$\lambda_1(O^{11}) = \lambda_2(O^{11}) = \gamma(1-m),\ \lambda_3(O^{11}) = -(1-l),\ \lambda_4(O^{11}) = -(1-l) - 2\alpha;$$
$$\lambda_{1,2}(AA) = \frac{mn \pm \sqrt{((mn)^2 - 4m\gamma v^*)}}{2},\ \lambda_{3,4}(AA) = \frac{mn - 2\alpha \pm \sqrt{((mn-2\alpha)^2 - 4m\gamma v^*)}}{2}$$

*where* $n = 1 + l - 2\alpha$.

Applying Statements 2, 3 and 4 of the Lemma to system (1s) we have the characteristic polynomial $\Psi(\lambda)$ in the form:

$\Psi(\lambda) \equiv \phi_4(\lambda)$ where $\phi_4(\lambda)$ is given by formula (9) at point $C(u_1,0,u_2,0)$,

$\Psi(\lambda) \equiv \varphi_4(\lambda)$ where $\varphi_4(\lambda)$ is given by formula (10) at point $B(y,0,m,z)$.

Expanding $\Psi(\lambda)$ in series on $\alpha$ and substituting the expansions of the non-symmetric equilibria coordinates, given in Theorem 2.1, we find eigenvalues of points we're considering.

***Corollary3.1.*** *Eigenvalues of System (1s) at the non-symmetric two-dimensional equilibria C are contained in Table 1. Eigenvalues of System (1s) at the non-symmetric three-dimensional equilibria B are contained in Table* 2.

Combining the previous results proves the next statement.

**Theorem 3.1.**

*For any parameter values belonging to the domain* ***M,*** *equilibria of System (1) have the following properties:*

*1)* $O$ *is a local asymptotically stable node,*

*2)* $O^{ll}, O^{11}$ *are unstable saddle-nodes in four-dimensional space;*

*3)* $C_{0l}, C_{01}, C_{l0}, C_{l1}, C_{10}, C_{1l}$ *are saddle-nodes that can have unstable as well as stable manifolds in their neighborhood* (see Table 1)*;*

*4)* $B(m, y, x, 0)$ *and* $B(x, 0, m, y)$ *-equilibria of the model (three pairs,* $B_1^0$, $B_1^l$, $B_1^1$ and $B_2^0$, $B_2^l$, $B_2^1$, *inside parameter domain bounded by* ***SB*** *and only one pair* $B_1$, $B_2$ *outside this domain) can be either a saddle-focus or a saddle-node; they are unstable in four-dimensional space.*

*5)* $AA(m, v^*, m, v^*)$ *where* $v^* = (m-l)(1-m)$ *is a stable spiral for* $(l+1)/2 < m < 1$, $\alpha > m(1+l-2m)/2$*, and an unstable spiral for* $0 < l < m$, $\alpha < m(1+l-2m)/2$.

## 4. Dynamics of "two/three-dimension subsystems" of Model (1s)

### 4.1. *Dynamics of "two/three-dimension subsystems" of Model (2)*

Now we are ready to describe the dynamics of 2D-system (3) and 3D-system (4) whose biological sense was discussed in the Introduction. The dynamics of model (2) are described mainly in the Introduction. The phase-parameter diagram, which was given in Fig. 2, has the following biological interpretation (for any fixed $\gamma > 0$). In domain 1 (*m*>1) predators go to extinction with any initial density, whereas the prey go extinct or to the steady state $u = 1$ depending on whether $u(t=0) < l$ or $u(t=0) > l$. In Domains 2 and 3 both prey and predator either coexist in steady state oscillations or go to extinction depending on initial densities. Finally, in Domains 4 and 5 both populations go to extinction for any initial data in spite of the fact that the model has positive non-trivial equilibrium in Domain 4.

System (3) is the simplest form displaying the model dynamics of two Allee-type prey populations interacting diffusely. The behavior of this system is defined completely by its equilibrium points. It has from three (for big $\alpha$) up to nine (for very small $\alpha$) equilibria. For any parameters this system has two stable nodes, $O_{00}(0,0)$ and $O_{11}(1,1)$, an unstable node $O_{ll}(l,l)$, and other equilibria that are saddles when they exist. The parameter-phase portrait of the model is given in Fig.3.

System (4) can be considered as a one-dimensional updating of System (3) in the case when a prey population has predators that control prey density. System (4) can be also considered as a one-dimensional updating of System (2) for the case when a prey population can be replenished with diffusely migrated preys. The bifurcation analysis of Model (4) provided by [6] and [7] have shown that the behavior of System (4) is much more complicated and diverse than the behavior of systems (3) and (2). Here we describe some main properties of these dynamics.

For $\alpha = 0$ and $\gamma = 1$ values $(l^*, m^*)$ are chosen in domain **4** of the parameter portrait of system (2) (see Fig.2) and give trajectories that tend to the origin from almost all initial data. (This case is interpreted as the community going to extinction). Let $B$ be an unstable non-trivial equilibrium of system (4). With increasing of parameter $\alpha$ the following behaviors were observed (see Fig.5).

*Statements* [6], [7].

*For any fixed* $0 < m < 1$ *and arbitrary* $0 < l < m$ *there exist parameters* $0 < \alpha^* < \alpha^{**} < \alpha^{***} < \alpha^{****}$ *dependent on* $l, m$ *such that*

*1) three-dimensional system (4) has no limit cycle for* $0 < \alpha < \alpha^* << 1$*;*

*2) for* $\alpha = \alpha^*$ *separatrices of the two two-dimensional equilibria* $C_{l0}, C_{10}$ $(C_{0l}, C_{01})$ *compose heteroclinics;*

*3) for* $\alpha^* < \alpha < \alpha^{**}$ *system (4) has a stable limit cycle* $c_3$*;*

*4) for parameter values belonging to the surface* $H^+{}_3 : \alpha^{**} \approx \dfrac{m(1+l-2m)}{2-m}$ *(see Table 2) equilibrium* $B$ *gains stability in a supercritical Andronov-Hopf bifurcation,*

*5) for* $\alpha^{**} < \alpha < \alpha^{***}$ *equilibrium* $B$ *is stable (see Fig.5);*

*6) for* $\alpha > \alpha^{****} > \alpha^{***}$ *the system has a unique stable manifold and equilibrium point at the origin.*

**5. Dynamics of System (1s)**

**5.1. *Change of stability of equilibrium*** $AA$

The (*l*,*m*)-parameter portrait of the local system (2), given in Fig.2, contains the boundary curve $\boldsymbol{H_1}$*:* $m = \frac{l+1}{2}$ corresponding to an supercritical Andronov-Hopf bifurcation of equilibrium $A(m, v^*)$ where $v^* = (m-l)(1-m)$. The eigenvalues of $A$ are: $\lambda_{1,2}(A) = (mn \pm \sqrt{((mn)^2 - 4m\gamma v^*)})/2$ where $n = 1 + l - 2m$. The equation of $\boldsymbol{H_1}$ is defined by the vanishing of the real parts of $l_{1,\,2}$.

An "image" of *A* in the frame of the bilocal system (1) is a point $AA(m, v^*, m, v^*)$ possessing two pair of eigenvalues:

$\lambda_{1,2}(AA) = (mn \pm \sqrt{((mn)^2 - 4m\gamma v^*)})/2,\ \ \lambda_{3,4}(AA) = (mn - 2\alpha) \pm \sqrt{((mn-2\alpha)^2 - 4m\gamma v^*)})/2$.

The vanishing of the real parts of $\lambda_{1,2}$ and $\lambda_{3,4}$ defines the "*neutrality*" surface ***H*** (see Fig. 6) consisting of two branches. There is the "old" branch

$$\boldsymbol{H^1}\{\gamma, \alpha, m, l : m = \frac{l+1}{2}\},$$

and the "new" branch (a parabola for any fixed $\gamma, l$):

$$\boldsymbol{H^2}\{\gamma, \alpha, m, l : \alpha = \frac{m(l+1-2m)}{2}\}.$$

***Proposition 5.1.*** *Let* $0 < l < m$. *Equilibrium* $AA(m, v^*, m, v^*)$ *of System (1s) located in the first orthant of the* $(u_1, v_1, u_2, v_2)$*-phase space, changes stability in supercritical Andronov-Hopf bifurcations with* *parameter values belonging to surface* $\boldsymbol{H^1}$ *for any* $\gamma > 0$, *and with parameter values belonging to surface* $\boldsymbol{H^2}$ *for any* $0 < \gamma \le 1$ .

*Proof.* We have calculated the first Lyapunov quantities $\boldsymbol{L_1}$ at $\boldsymbol{H^1}$ and $\boldsymbol{H^2}$ (see [3], [19]). Applying the procedure contained in [15] we obtained the following results.

Let $v_0 = (m-l)(1-m)$ is $v_1 = v_2$-coordinate of equilibrium *AA*. Then

$$L_1(H^1) = -c_1 \frac{1+l}{2\gamma v_0},$$

and

$$L_1(H^2) = c_2(-8\alpha^2(4+4l-\gamma-6m) - 9(1+l)\gamma m v_0 - 2\alpha((8+8l^2+2\gamma-\gamma^2+2l(8+\gamma))m - 6(8+8l+\gamma)m^2 + 72m^3 + 2\gamma m v_0)),$$

*where* $\alpha = m(l+1-2m)/2.$

Here $c_1$, $c_2$ are positive constants.

It is evident that $L_1(H^1) < 0$. Analysis of formula $L_1(H^2)$ has been done numerically. We found $L_1(H^2) < 0$ and have also verified the results using the LOCBIF–packages [17]. This completes the proof.

We note that our computer experiments with system (1) revealed that the stable cycle, which appears when crossing the boundary $\boldsymbol{H^2}$, disappears with parameter values very close to those in $\boldsymbol{H^2}$, annihilating with an unstable limit cycle. Because of this we did not mention the domain of its existence in the parameter portrait of system (1s), see Fig.7, and denote $\boldsymbol{H^2}$ by dotted line there. The stable cycle, which appears crossing the boundary $\boldsymbol{H^1}$ exists a wide parameter domain; we call this cycle $\boldsymbol{c_u}$.

5.2. **4D-oscillations in System (1s)**

Computational analysis of the system (1s) reveals the complicated structure of its dynamics, which essentially depends on parameter values of the system as well as on initial values of the variables. Of course, due to biological interpretations of the model (1s) our main interest is in the *stable* modes, which can be observed with variations of initial data. We show that even for fixed parameters the system can demonstrate wide range multistability. In addition, the existence of migration in our model increases diversity of stable modes thereby increasing the sustainability of the model community.

The schematic parameter portrait of the system for some typical 4D-rearrangements of model behaviors in a phase neighborhood of point $AA$ is presented in Fig. 6. The portrait represents an $(\alpha, m)-$cut of the four-dimensional $(\alpha, \gamma, m, l)-$ parameter space for fixed values $\gamma = 1$ and $0 < l = .1 < m < 1$, such that $(l, m)$ belongs to Domain **4** or adjacent to the boundary ***L*** Domain **3** in Fig.2. The portrait was obtained by analytical and computer methods of the bifurcation theory with the use of packages TRAX [21] and LOCBIF [17]. The analytically obtained curves $\boldsymbol{H}^1$, $\boldsymbol{H}^2$ of Andronov-Hopf bifurcations were described above. The 3D-rearrangements of model behavior (see Fig.5) are not presented in this portrait and will be discussed later.

In Domain **I** of the parameter portrait the model has stable four-dimension equilibrium *AA*.

In Domain **II** of the parameter portrait the model has stable oscillations corresponding to the stable four-dimensional limit cycle "$\boldsymbol{c_u}$" in $(u_1, v_1, u_2, v_2)$- phase space. This cycle appears when we cross boundary $\boldsymbol{H}^1$ from bigger to smaller $m$ (from right to left in Fig.6) and exists for any value of parameter $\alpha \geq 0$ (see Fig. 7). Its "pre-image" in system (2) also appeared in a supercritical Andronov-Hopf bifurcation after crossing boundary ***H*** (see Fig.2). Recall that ***H*** and $\boldsymbol{H}^1$ have the same equation, and equilibrium $AA$ loses stability in the two-dimensional eigen-space at boundary $\boldsymbol{H}^1$, whereas the other two-dimensional eigen-space corresponds to trajectories tending to *AA* with $t \to \infty$. Computer analysis of cycle $\boldsymbol{c_u}$ (see Fig. 7) has shown that $\boldsymbol{c_u}$ disappears at heteroclinics composed by separatrices of the trivial symmetric equilibria $O^{ll}, O^{11}$, analogous to its pre-image which has disappeared at heteroclinics of $O_l, O_1$. The parametric boundary corresponding to the mentioned heteroclinics is denoted $\boldsymbol{L}^1$ in the portrait and schematically presented as a straight line. With fixed value of $l$ lines $\boldsymbol{L}^1$ in Fig.6 and $\boldsymbol{L}$ in Fig.2 has a common point for $\alpha = 0$: $(l = .1, m \approx .446)$.

In Domain **III** the model has only one 4D-attractor –equilibrium at the origin.

Model (1s) demonstrates its most interesting 4D-behavior in Domains **IV** and **V** bounded by curve $\boldsymbol{F}$ at the portrait (see Fig.4). Our computer analysis has revealed that for very small $\alpha$ (close to the boundary $\boldsymbol{SC_1}$ in Fig.3) and certain *m,* separatrices of trivial two-dimensional equilibria $C_{10}, C_{01}$ gives rise to four-dimensional limit cycle $\boldsymbol{c_4}$ (see Fig.8a, where $\alpha = .001$ and both cycles $\boldsymbol{c_4}$ and $\boldsymbol{c_u}$ are presented). For fixed values of $l$ increasing $\alpha$ results in $\boldsymbol{c_4}$ undergoing period-doubling, whereby it loses its stability transforming to a torus, etc., depending on the value of $m$ for which we consider this cycle. Fig.-s 9a and 9b, where $m = .45$ , $\alpha = .0256$ and $\alpha = .0235$ respectively, display the stages of period-doubling of $\boldsymbol{c_4}$ in Domain **IV**. In Fig.9b one can observe period-doubling in more detail; it is clear from this picture that domains of attraction of $\boldsymbol{c_u}$ and $\boldsymbol{c_4}$ are divided by an unstable manifold similar to a limit cycle. Note that $\boldsymbol{c_4}$ is destroyed at domain **IV** under some period-doubling, presumably, Feigenbaum doubling that “originates” a cycle of period 3 which leads to the appearance of weak chaos dynamics [27], [22]. In considering the model with $l = .1$ the destruction of cycle $\boldsymbol{c_4}$ was observed at the boundary $\boldsymbol{F}$ for $m = .48, \alpha \approx .0183$ **(**see Fig. 10a) as well as for $m = .446, \alpha \approx .0256$ **(**see Fig. 10b); note, that these pictures contain both cycles $\boldsymbol{c_4}$ and $\boldsymbol{c_u}$.

In Domain **V** increasing $\alpha$ results in cycle $\boldsymbol{c_4}$ transforming to a torus (when $m$ decreases, see Fig-s. 11, 12). This torus possessing irrational rotating number (see Fig.12, *m*=.325 with $\alpha$=.0332) is destroyed by the appearance of a rational rotation number (see Fig. 12b, $m$=.325 with $\alpha$=.0333) and chaotic dynamics result (see [14] for theory of general torus destruction).

In Fig.13 we show portraits of the system for $m=.315, \alpha=.027$; this parameter point is placed close to the boundary ***F*** of Domains **V** and **III**. The portraits display two stable manifolds: two-leaf torus $\boldsymbol{c_4}$ and three-dimensional cycle $\boldsymbol{c_3}$***;*** initial points in their basins are distinguished only in their $u$-coordinate, $u$=.01 for the former and $u$=0 for the latter.

Recall that for small enough $\alpha$ and $m < \frac{1+l}{2}$ system (1s) has two stable "three-dimensional" limit cycles $\boldsymbol{c_3^1}$ and $\boldsymbol{c_3^2}$ at spaces $(u_1, v_1, u_2, 0)$ and $(u_1, 0, u_2, v_2)$ (see Fig. 5).These cycles appear from heteroclinics composed by two-dimensional trivial equilibria $C_{l0}, C_{10}$ and $C_{0l}, C_{01}$ (at value $\alpha = \alpha^*$ in Fig. 5) and disappear at supercritical Andronov-Hopf bifurcation with parameter value $\alpha = \alpha^{**}$ in Fig. 4. Note that cycles $\boldsymbol{c_3^1}$ and $\boldsymbol{c_3^2}$ continue to exist in subspaces $(u_1, v_1, u_2, 0)$ and $(u_1, 0, u_2, v_2)$ of space $(u_1, v_1, u_2, v_2)$ and are stable if the initial values belong to these subspaces (see, for example, Fig.-s 10b and 13). Thus, for a wide range of parameters system (1s) simultaneously has 4D- and 3D-attractors, limit cycles. The origin equilibrium is stable for any parameter values, though its basin of attraction can vary. Thus, the model community can tend to one of these depending on initial values.

### 6. Discussion. Biological Interpretations

It is useful to compare dynamics of the three models described by the two-dimensional System (2), three-dimensional System (4), and four-dimension System (1s).

The model of the local (one patch) community (2) predicts four different regimes of dynamical behavior: getting predators to extinction with any initial density because the death rate of predators is too large; possibility of predator-prey coexistence at steady state or stable oscillations for a wide range of initial densities; getting preys and predators to extinction because predators over-regulate prey density.

The model of two-patch community (4) consisting of preys and predator-preys systems can have up to thirteen equilibrium points, both trivial (with at least one zero coordinate) and non-trivial. For suitable dispersal prey rates between both patches, the predator can control the numbers of both prey populations in a stable stationary equilibrium or in oscillatory regimes. It was revealed also that "trivial" two-dimensional equilibrium points in the frame of Model (4) "originate" oscillations in the community after their separatrices compose heteroclinics. The most "interesting" and diverse behaviors of this model are observed close to boundaries of population coexistence. Regimes of bi-stability or tri-stability were found: depending on initial data the model community can coexist at a non-trivial stable equilibrium, at stable oscillations, or go extinct. It is important to note that both systems in the model community can survive for parameter values for which any one of them would go to extinction for all initial densities in isolation.

The Model (1s) describing dynamics of the two-patch predator-prey system supports more diversity of stable modes when compared with Model (4). In fact, it can have up to sixteen equilibria and up to seven stable limit cycles in the 1-st orthant. We have analyzed Model (1s) mainly for small values of the dispersal parameter $\alpha$ with fixed values for the Malthusian parameter $\gamma \approx 1$. The main attention was on the investigation of the 4D-modes. Results of the analyses are collected in the schematic parameter portrait in Fig.6 as an $(\alpha, m)-$cut of the four-dimensional $(\alpha, \gamma, m, l)-$parameter space obtained by fixing $\gamma = 1$ and $0 < l = .1 < m < 1$. We follow the modes of the stable community dynamics by changing the parameters $m$ and $\alpha$. The parameter portrait of the system is divided into 5 domains comprising different 4D-phase behaviors of the model. We found the boundaries $\boldsymbol{H^1}$, $\boldsymbol{H^2}$ corresponding to the change of stability of non-trivial 4D-equilibria in a supercritical Andronov-Hopf bifurcation. The boundaries $\boldsymbol{L}^1$ and $\boldsymbol{F}$ in the portrait correspond to non-local bifurcations in the model, “heteroclinics” of equilibrium points, and limit cycles. We showed, for example, that for parameter values belonging to domain **IV** the model community can support up to **four stable nontrivial modes** under suitable initial values: oscillations corresponding to two 4D limit cycles as well as those corresponding to two 3D limit cycles were found. The size of the basins of attraction of these cycles depend on parameters, we observed and analyzed the changing of form and the range of the oscillations when $\alpha$ increased. We observed period doubling leading to aperiodic, possibly chaotic oscillations similar to those described in [27] and [22]. In addition, we observed the formation and destruction of a torus as the rotation number moves from irrational to rational (a known route to chaos—see [14]). Note, that decreases in the parameter *m,* which lead to the extinction of the predator-prey system in one patch, do not necessarily lead to the “immediate” extinction of the community for suitable values of $\alpha$. In fact, *m* and $\alpha$ combinations result in communities that persist periodically, aperiodically, or chaotically oscillating.

The two limit cycles $\boldsymbol{c_u}$ and $\boldsymbol{c_4}$ can be biologically interpreted as follows: 4D-“oscillations $\boldsymbol{c_u}$” are originated from 2D- oscillations that live in one-patch Models (2), persist for any $\alpha$ in Domain II (in the frame of (1s)). 4D-“oscillations $\boldsymbol{c_4}$” arise due to the structure of Model (1s) for very small $\alpha$ and reflect the ability of populations to coexist in a rather extreme oscillatory regime, under weak dispersal, see the limit cycle in Fig. 8b. As the parameter $\alpha$ increases or the parameter $m$ decreases these regular oscillations move into a period-doubling regime or in a torus – a typical route to destruction in the chaotic regime. In other words, excessive dispersal or low mortality predator death (given predators the opportunity to eat all the prey) leads to the system extinction in a rather exotic way.

3D-oscillations, which are one of the stable modes of a system of one predator- two preys (see Model (4)), coexist simultaneously with 4D-oscillations; small changes in the predators' initial densities can shift the dynamics from 3D-stable oscillations to 4-D stable oscillations and vice-versa. 3-D oscillations are more robust, that is, the 4D-oscillations may die while the 3-D survive.

These model outcomes also hold, supported by numerical analysis, in the non-symmetric version of the model. We also considered the case where the Allee parameters $l_1$ and $l_2$ were not identical. There we observed the appearance of new stable modes together with the old one (see Fig.14a,b).

Further investigation can follow by considering a fully non-symmetrical system (1).

**Acknowledgment.**

The authors greatly appreciate to Dr. N. Davydova for her important contribution to this work being a student of F.B. The authors express gratitude for the assistance in computations provided by Dr. Gordillo and the help in finding the first Lyapunov values provided by Dr. Novozhilov. This project has been partially supported by grants from the National Science Foundation (NSF - Grant DMPS-0838704), the National Security Agency (NSA - Grant H98230-09-1-0104), the Alfred T. Sloan Foundation and the Office of the Provost of Arizona State University.

Fig.1. Schematic for predator-prey interactions. Note that although preys are allowed to migrate between patches, the predators are *not* allowed to migrate. In this paper, we consider the following situations: eq.(1) is $\alpha_i>0, \gamma_i>0$; eq.(2) is $\alpha_i=0, \gamma_i>0$; eq.(3) is $\alpha_i>0, \gamma_i=0$; eq.(4a) is $\alpha_i>0, \gamma_1>0, \gamma_2=0$; eq.(4b) is $\alpha_i>0, \gamma_1=0, \gamma_2>0$.

Fig.2. Schematically presented the parameter-phase portrait of model (2).

Fig.3. Parameter (a) and corresponding phase (c) portraits of non-symmetric equilibria $C(u_1^*, u_2^*)$ for model (3), non-symmetric trivial equilibria $C(u_1^*, 0, u_2^*)$ for system (4a), $C(u_1^*, u_2^*, 0)$ for system (4b) and $C(u_1^*, 0, u_2^*, 0)$ for system (1s) correspondingly, where $u_1^*, u_2^*$ are roots of (5). There are three pairs of these equilibria in Domain 1, one pair in Domain 2, no equilibria in Domain 3. Fig. (b) explains the notation of *C*-equilibria, $C_{0l}, C_{l0}, C_{01}, C_{10}, C_{l1}, C_{1l}$, as those arising from two of the equilibria $O_0, O_l, O_1$ of system (2) affected by parameter $\alpha$. $C_{00} \equiv O_{00}$, $C_{ll} \equiv O_{ll}$, $C_{11} \equiv O_{11}$ are also presented in the pictures; they exist for any $\alpha$.

Fig. 4. Parameter (a) and corresponding phase (b) portraits of non-trivial equilibria $B(y, m, z)$ for Model (4b), $B(m, z, y)$ for model (4a), trivial equilibria $B_1(y, 0, m, z)$ and $B_2(m, z, y, 0)$ for system (1s) correspondingly, where $y, z$ are roots of (7). System (1s) has three pairs of these equilibria in Domain 1 ($B_1^0$, $B_1^l$, $B_1^1, B_2^0$, $B_2^l$, $B_2^1$) and one pair in Domain 2 ($B_1, B_2$). The boundary between domains corresponds to the fold bifurcation in any points except upper point corresponding to the cusp bifurcation.

Fig.5. The main stable modes of dynamics of model (4) when $\alpha$ changes

Fig.6. (**a**) Schematically presented $(\alpha, m)-$cut of the $(\gamma = 1, l = .1, \alpha, m)-$ parameter portrait of 4D-stable modes of model (1). Domain **I** contains the stable non-trivial 4D-equilibrium, Domain **III** has no non-trivial 4D-attractors. (**b**) Phase portraits in domains **II**, **IV** and **V.**

Fig.7. Limit cycle $c_u$ appears in the subcritical Hopf bifurcation with $m = .55$ and disappears in heteroclinics with $m \cong .45$

Fig.8. In domain **IV** limit cycles $c_4$ and $c_u$ coexist for parameters $\alpha = .001$, $m$=.46, $\gamma = 1, l = .1$. Limit cycles $c_4$ and $c_u$ are shown in plane $(u_1, u_2)$ (**a**), and in plane $(u_1 + u_2, v_1 + v_2)$ (**b**)

Fig.9. Period-doubling of cycle $c_4$ in domain **IV;** (a): $\alpha = .0256$, (b): $\alpha = .0235$; here cycle $c_u$ is also presented.

Fig.10. Cycle $c_4$ destructs at the boundary ***F:***(a) for $m = .48, \alpha \approx .0183$, (b) for $m = .446, \alpha \approx .0256$

Fig.11. Changing cycle $c_4$ with decreasing of parameter $m$ in Domain **V.**

Fig.12. Irrational rotating number at $c_4$ in domain **V**

Fig.13.Stable"cycles" $c_4$ and $c_3$ coexists in domain **V**

Fig.14. Examples of dynamical modes arising in the model if parameters $l_1 \neq l_2$

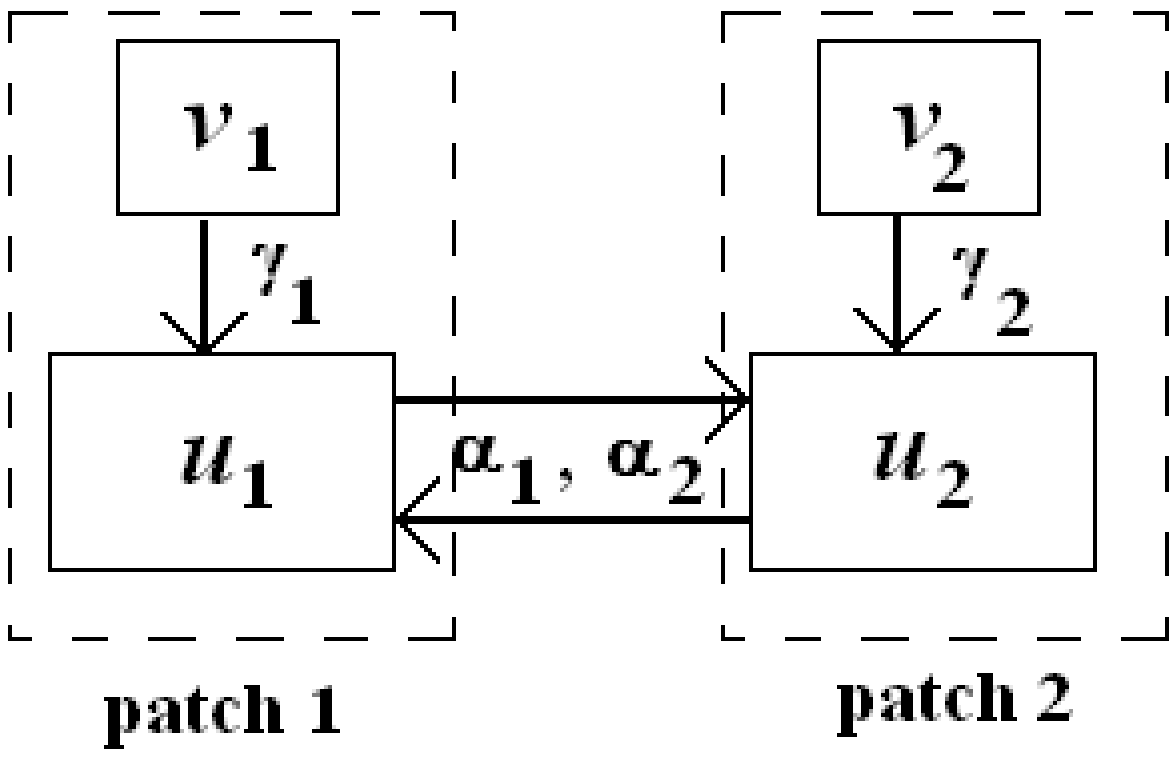

Fig.1. Schematic for predator-prey interactions. Note that although preys are allowed to migrate between patches, the predators are *not* allowed to migrate. In this paper, we consider the following situations: eq.(1) is $\alpha_i>0$, $\gamma_i>0$; eq.(2) is $\alpha_i=0$, $\gamma_i>0$; eq.(3) is $\alpha_i>0$, $\gamma_i=0$; eq.(4a) is $\alpha_i>0$, $\gamma_1>0$, $\gamma_2=0$; eq.(4b) is $\alpha_i>0$, $\gamma_1=0$, $\gamma_2>0$.

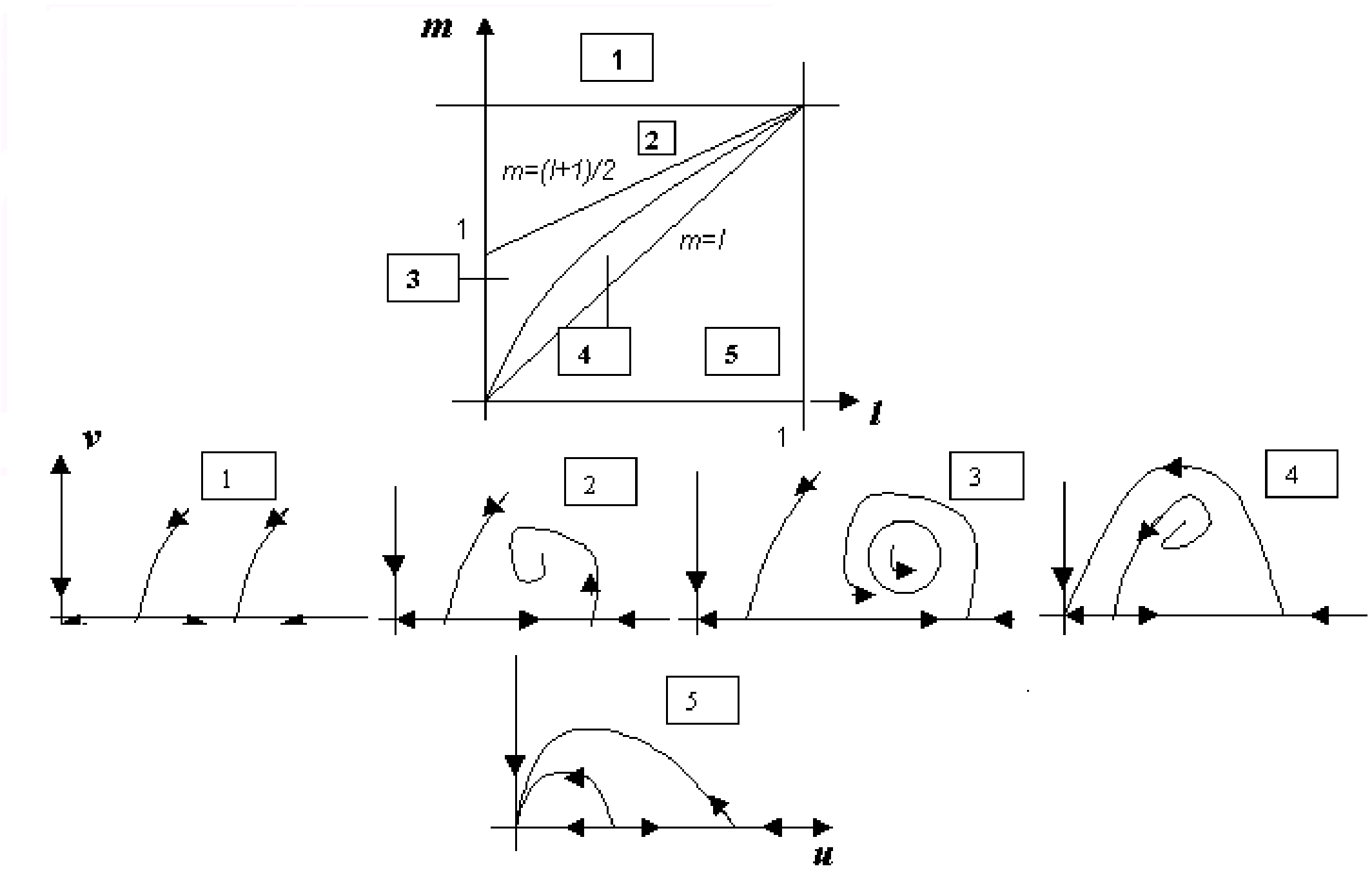

Fig.2. Schematically presented the parameter-phase portrait of model (2).

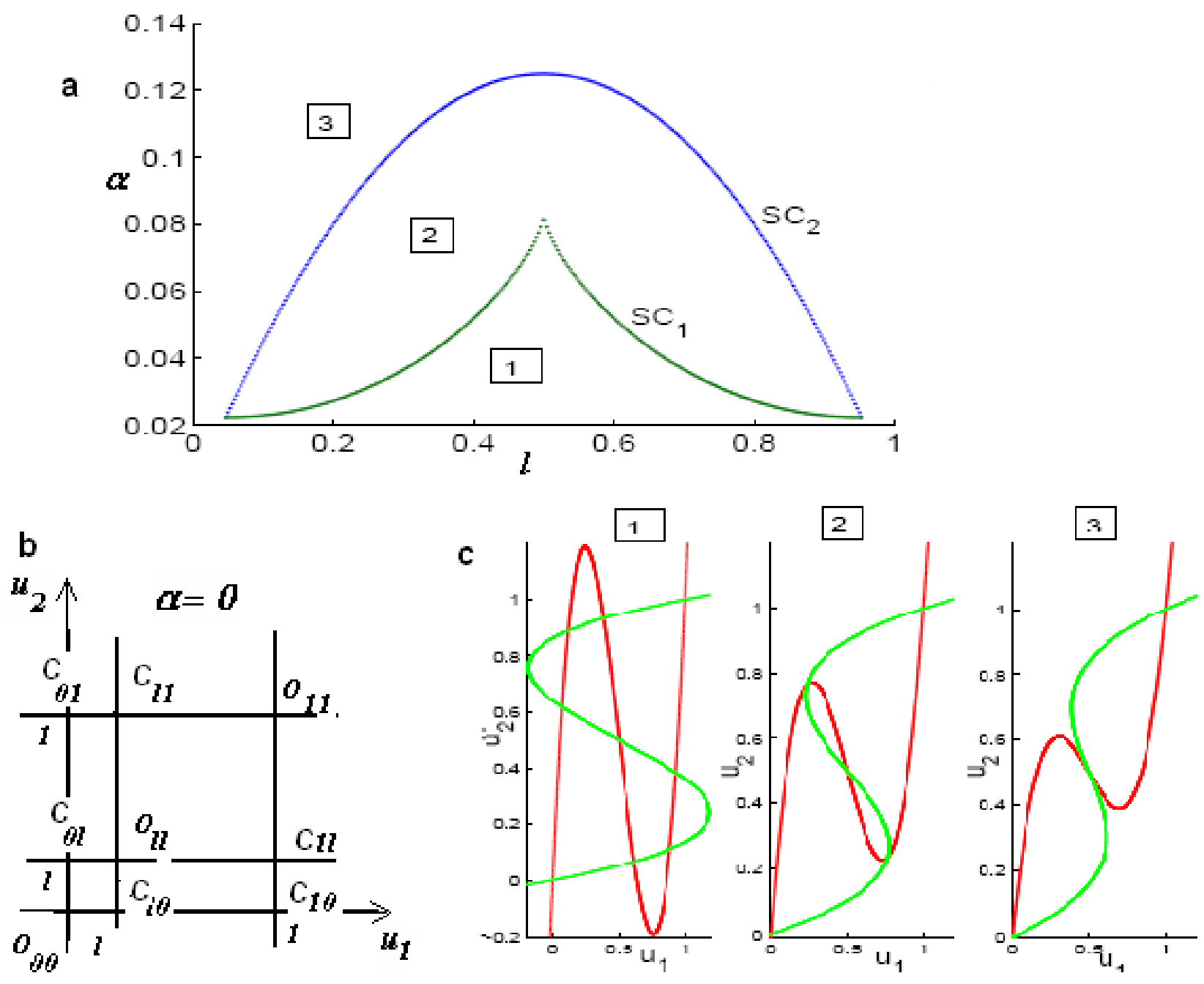

Fig.3. Parameter (a) and corresponding phase (c) portraits of non-symmetric equilibria $C(u_1^*, u_2^*)$ for model (3), non-symmetric trivial equilibria $C(u_1^*, 0, u_2^*)$ for system (4a), $C(u_1^*, u_2^*, 0)$ for system (4b) and $C(u_1^*, 0, u_2^*, 0)$ for system (1s) correspondingly, where $u_1^*, u_2^*$ are roots of (5). There are three pairs of these equilibria in Domain 1, one pair in Domain 2, no equilibria in Domain 3. Fig. (b) explains the notation of *C*-equilibria, $C_{0l}, C_{l0}, C_{01}, C_{10}, C_{l1}, C_{1l}$, as those arising from two of the equilibria $O_0, O_l, O_1$ of system (2) affected by parameter $\alpha$. $C_{00} \equiv O_{00}$, $C_{ll} \equiv O_{ll}$, $C_{11} \equiv O_{11}$ are also presented in the pictures; they exist for any $\alpha$.

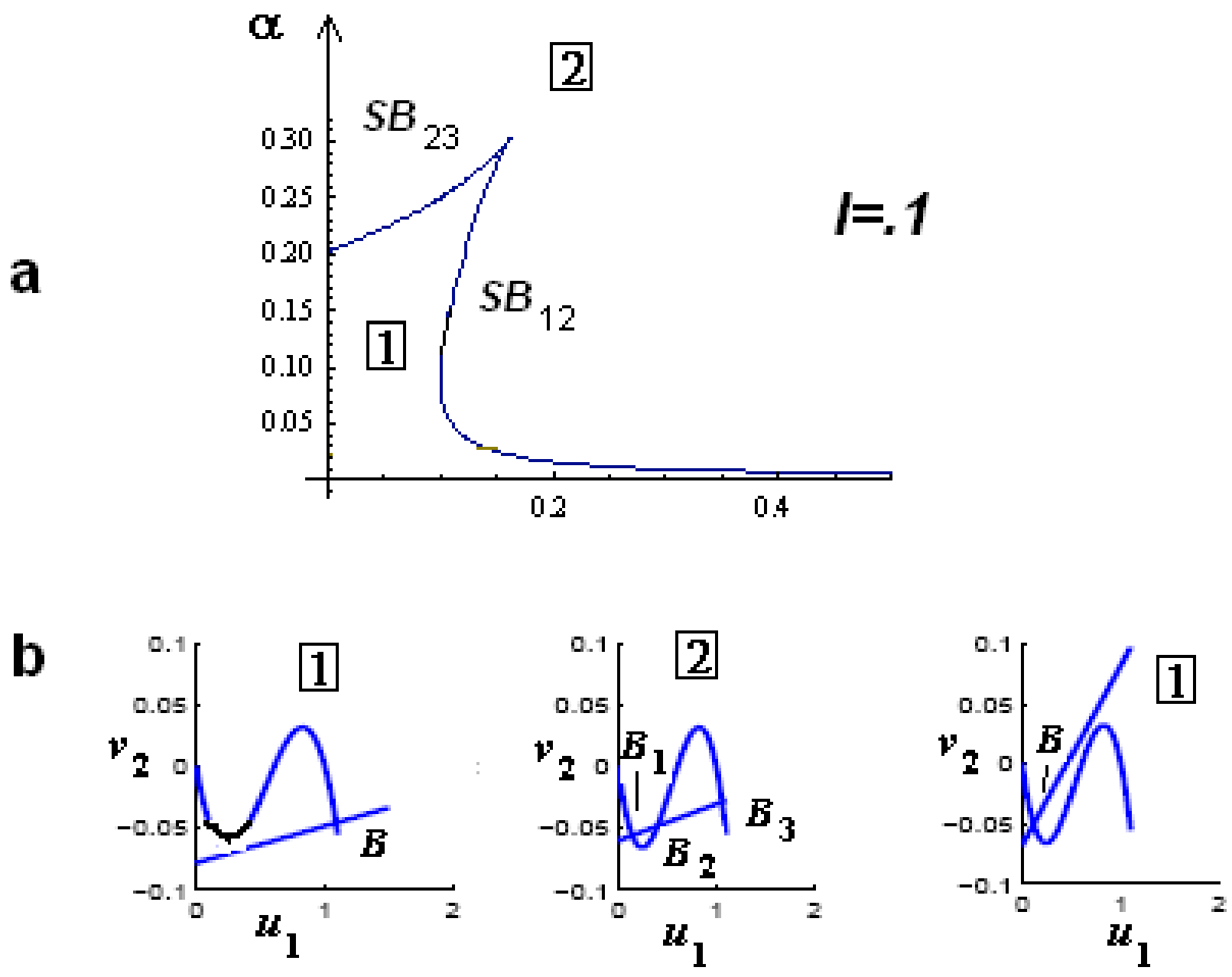

Fig. 4. Parameter (a) and corresponding phase (b) portraits of non-trivial equilibria $B(y,m,z)$ for Model (4b), $B(m,z,y)$ for model (4a), trivial equilibria $B_1(y,0,m,z)$ and $B_2(m,z,y,0)$ for system (1s) correspondingly, where $y,z$ are roots of (7). System (1s) has three pairs of these equilibria in Domain 1 ( $B_1^0$, $B_1^l$, $B_1^1$, $B_2^0$, $B_2^l$, $B_2^1$) and one pair in Domain 2 ( $B_1$ , $B_2$ ). The boundary between domains corresponds to the fold bifurcation in any points except upper point corresponding to the cusp bifurcation.

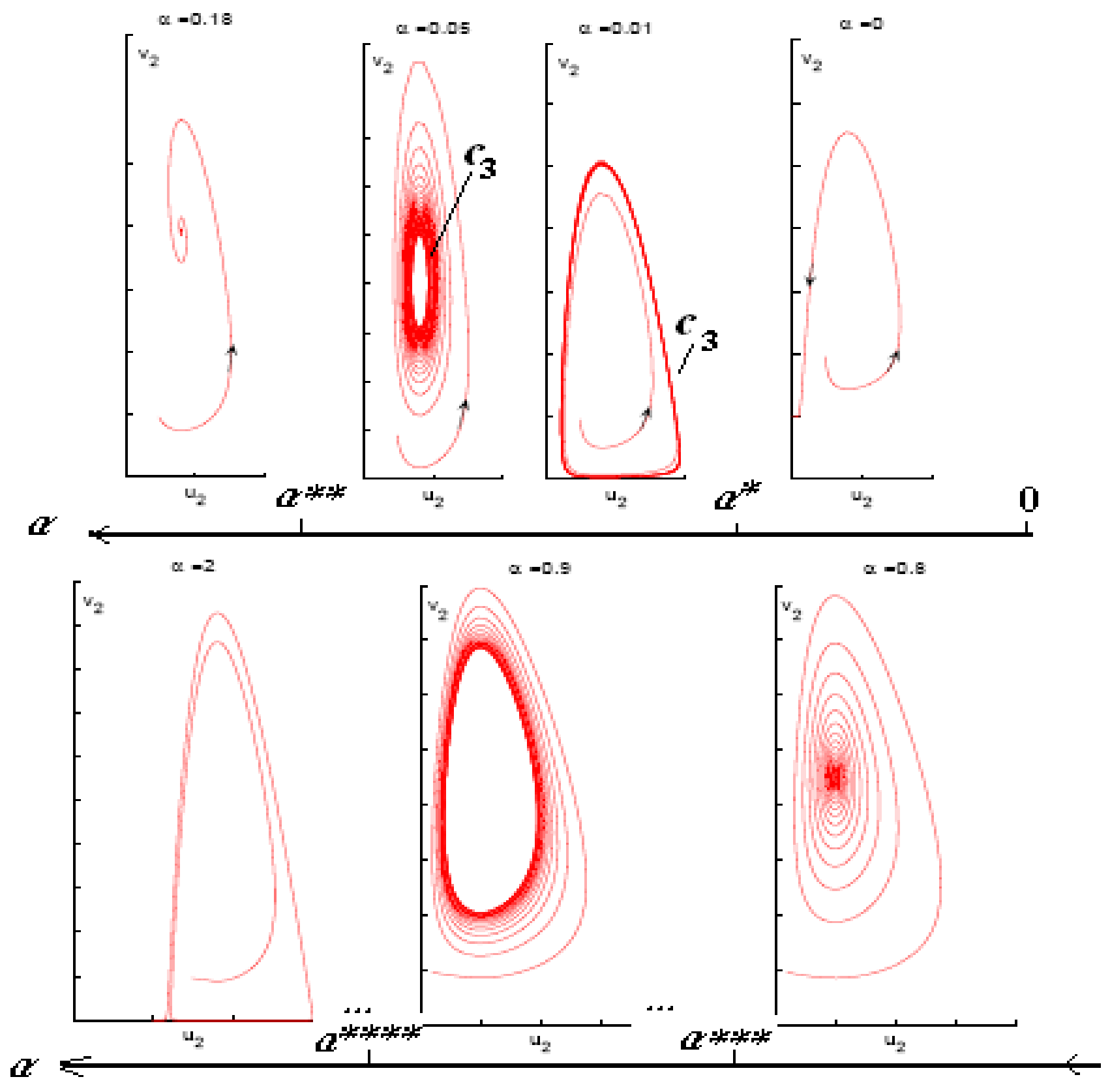

Fig.5. The main stable modes of dynamics of model (4) when $\alpha$ changes

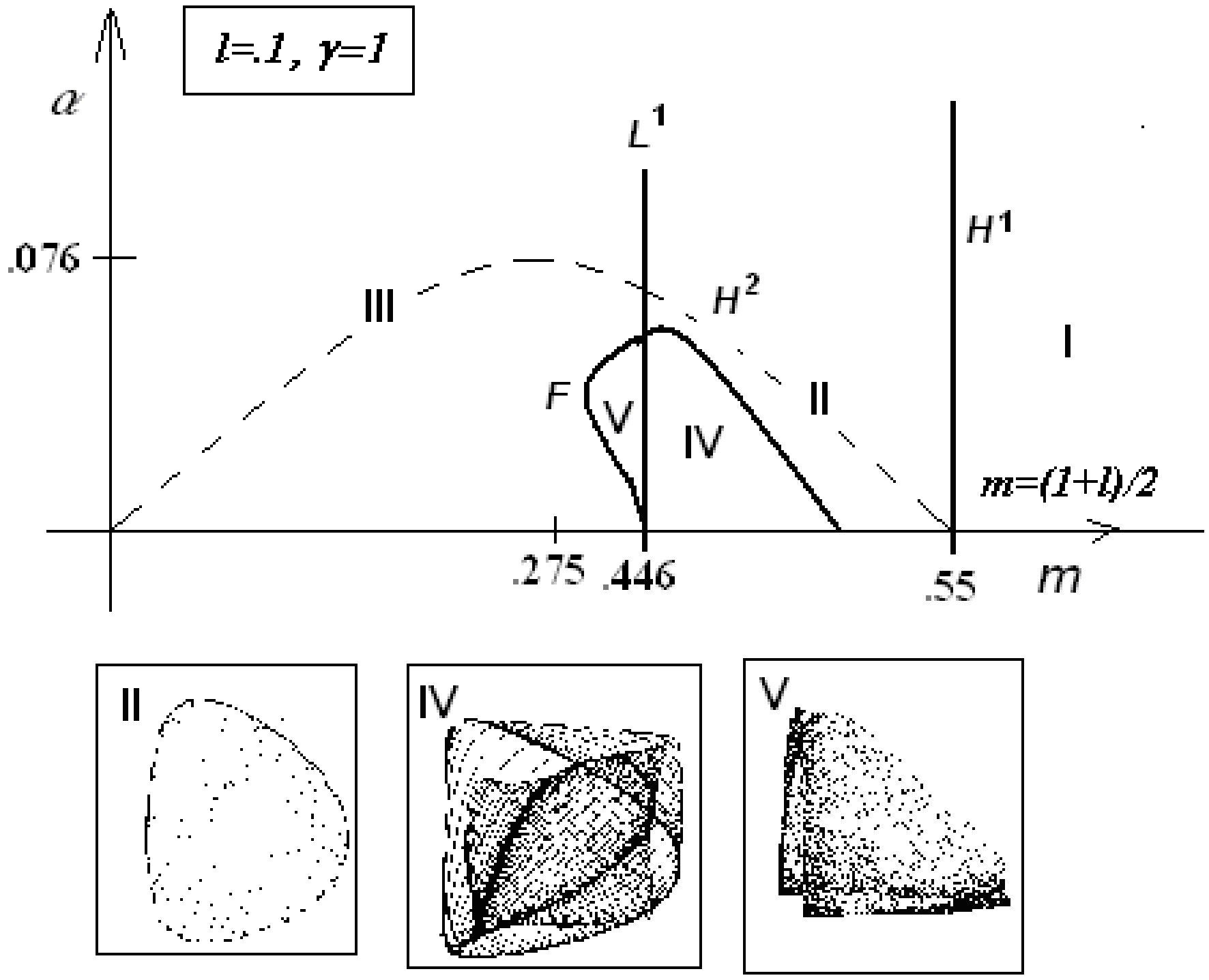

Fig.6. (**a**) Schematically presented $(\alpha, m)-$cut of the $(\gamma=1, l=.1, \alpha, m)-$ parameter portrait of 4D-stable modes of model (1). Domain **I** contains the stable non-trivial 4D-equilibrium, Domain **III** has no non-trivial 4D-attractors. (**b**) Phase portraits in domains **II**, **IV** and **V.**

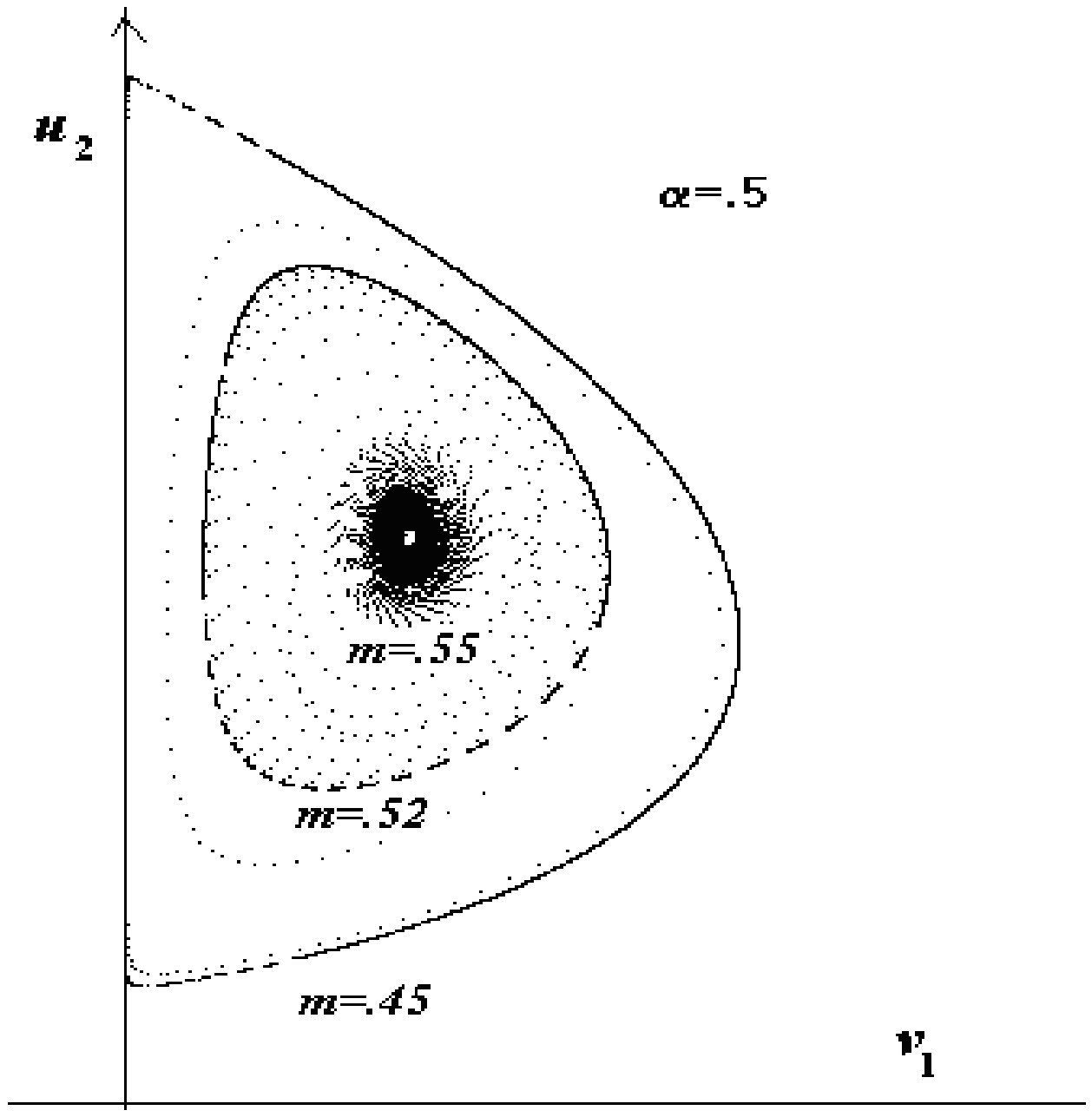

Fig.7. Limit cycle $\boldsymbol{c_u}$ appears in the subcritical Hopf bifurcation with $m = .55$ and disappears in heteroclinics with $m \cong .45$

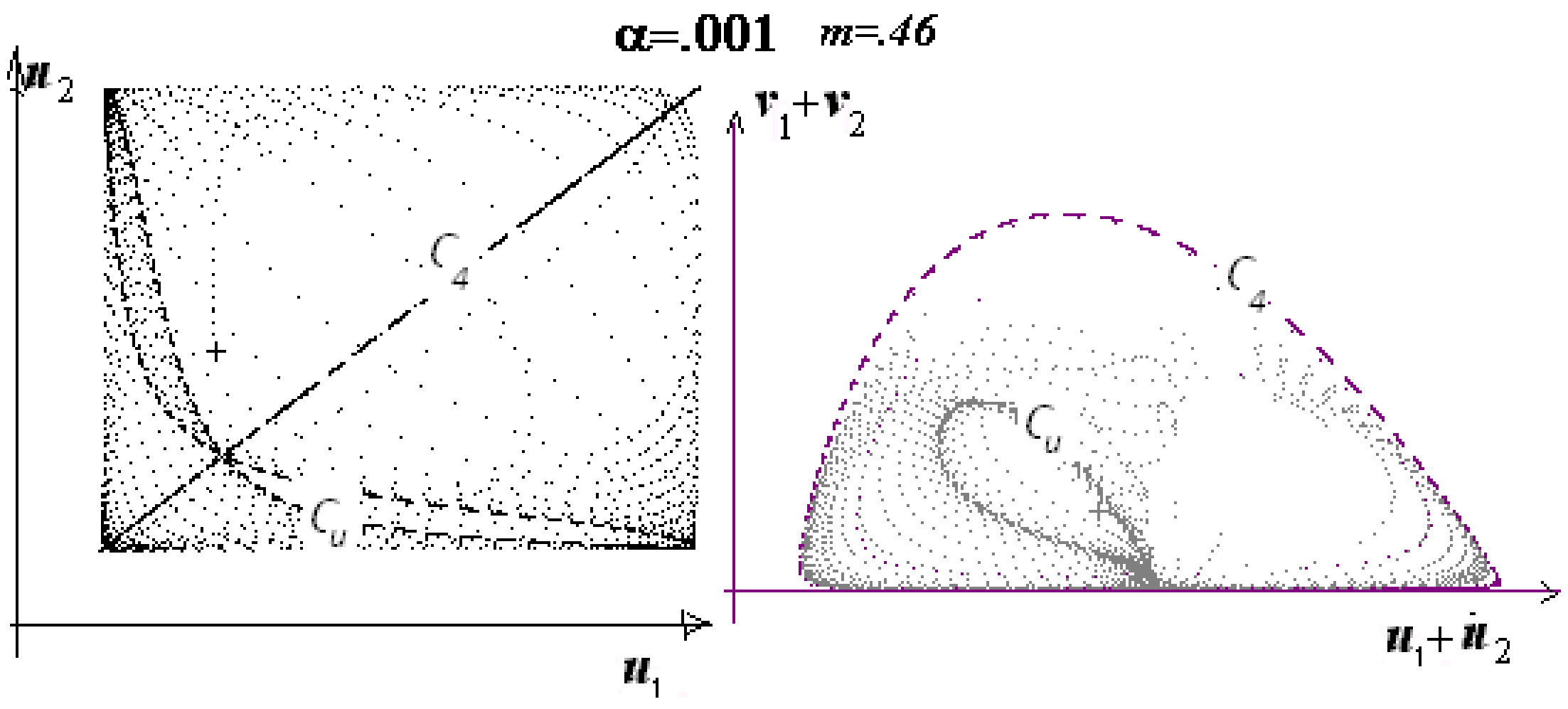

**a** **b**

Fig.8. In domain **IV** limit cycles $\boldsymbol{c_4}$ and $\boldsymbol{c_u}$ coexist for parameters $\alpha = .001$, $m$=.46, $\gamma = 1, l = .1$. Limit cycles $\boldsymbol{c_4}$ and $\boldsymbol{c_u}$ are shown in plane $(u_1, u_2)$ (**a**) , and in plane $(u_1 + u_2, v_1 + v_2)$ (**b**)

**a**

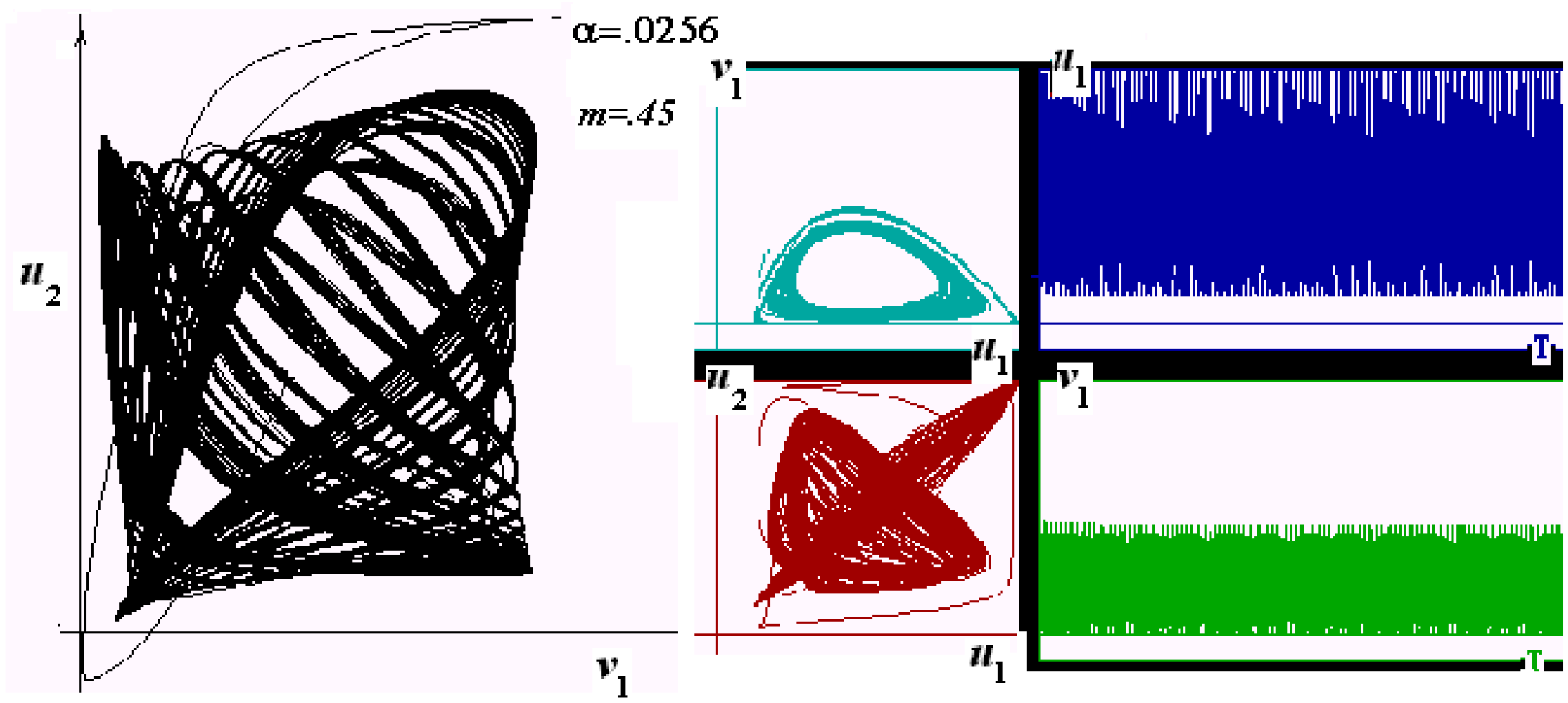

**b**

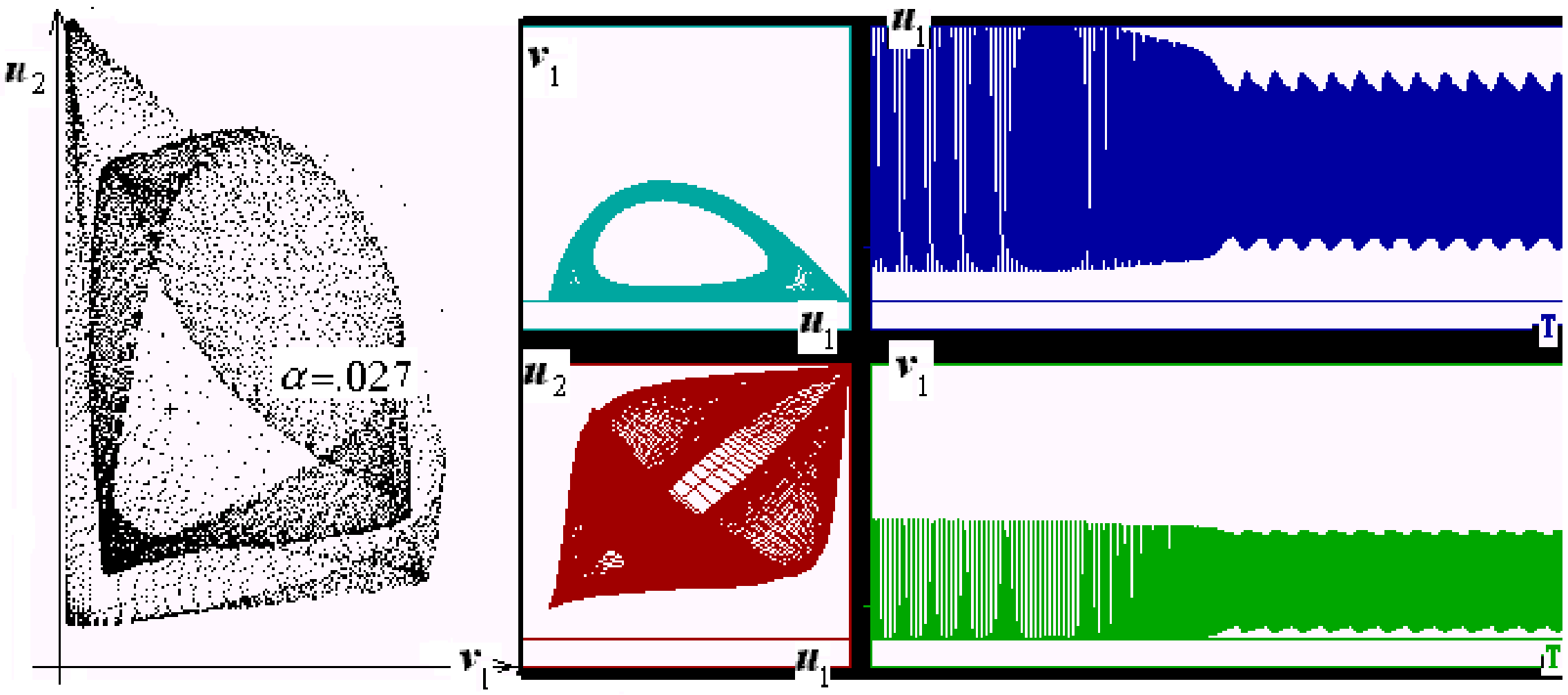

Fig.9. Period-doubling of cycle $c_4$ in domain **IV;** (a): $\alpha = .0256$, (b): $\alpha = .0235$; here cycle $c_u$ is also presented.

**a** **b**

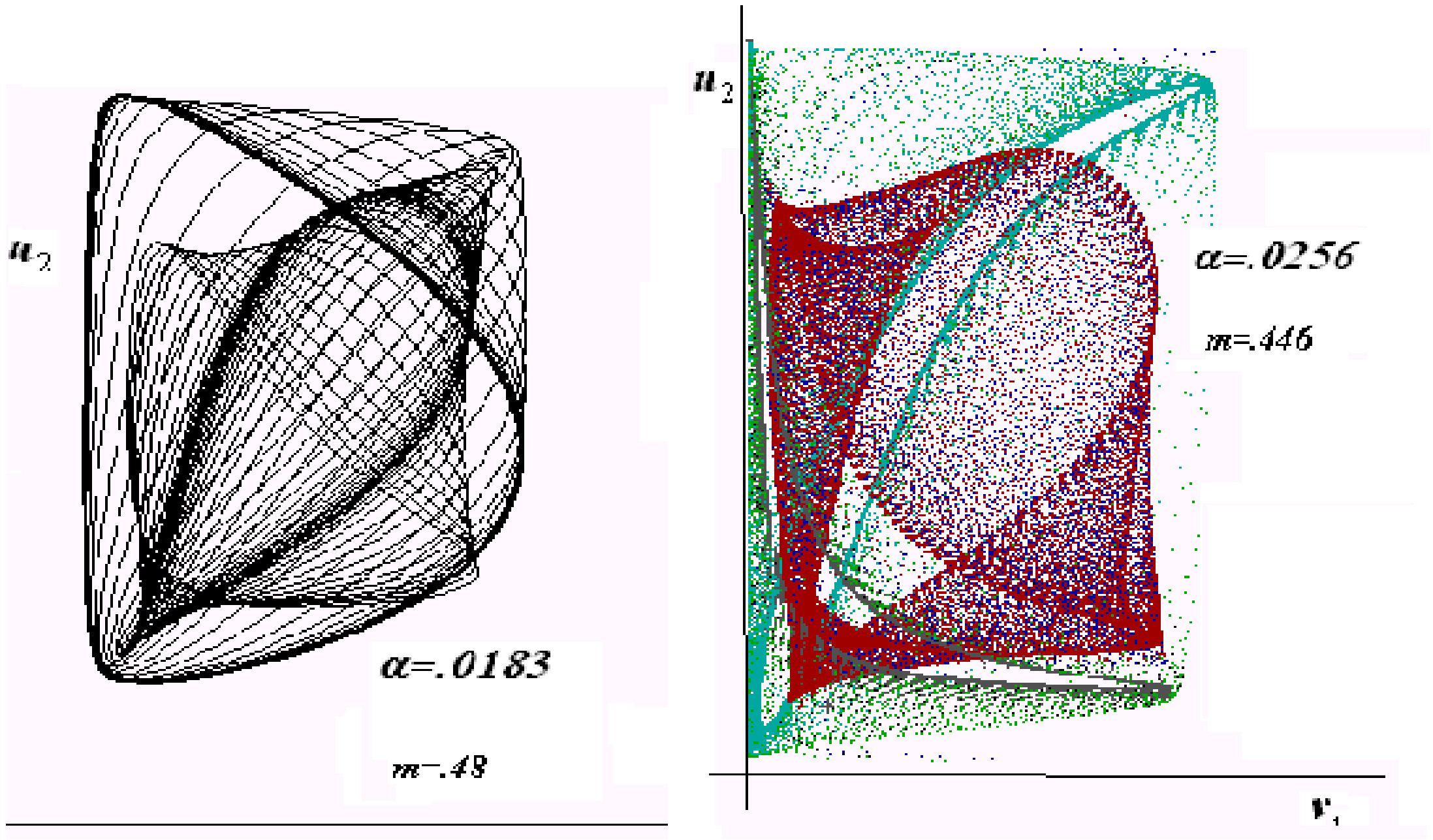

Fig.10. Cycle $c_4$ destructs at the boundary $F$:(a) for $m = .48, \alpha \approx .0183$, (b) for $m = .446, \alpha \approx .0256$

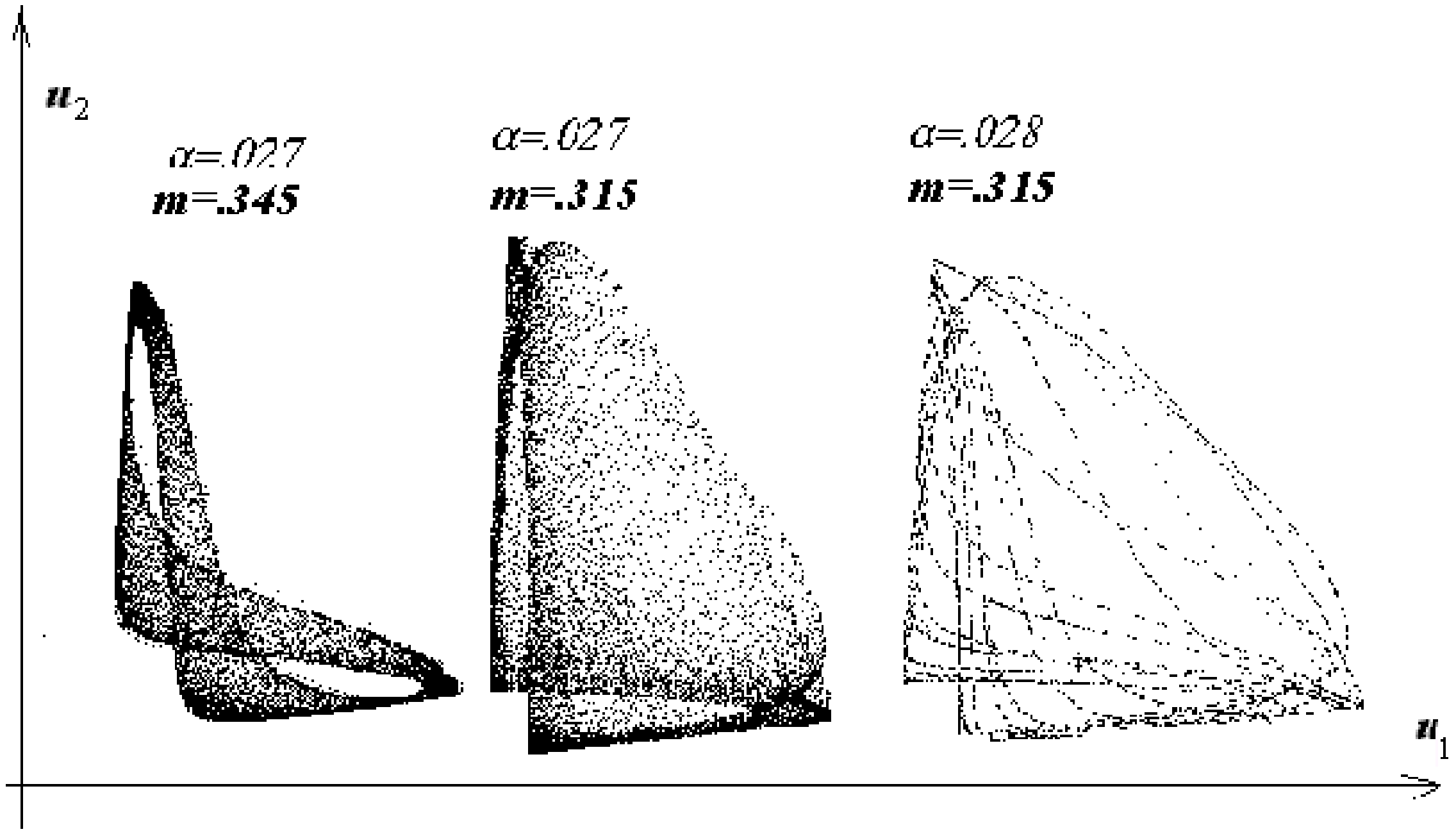

Fig.11. Changing cycle $c_4$ with decreasing of parameter $m$ in Domain **V.**

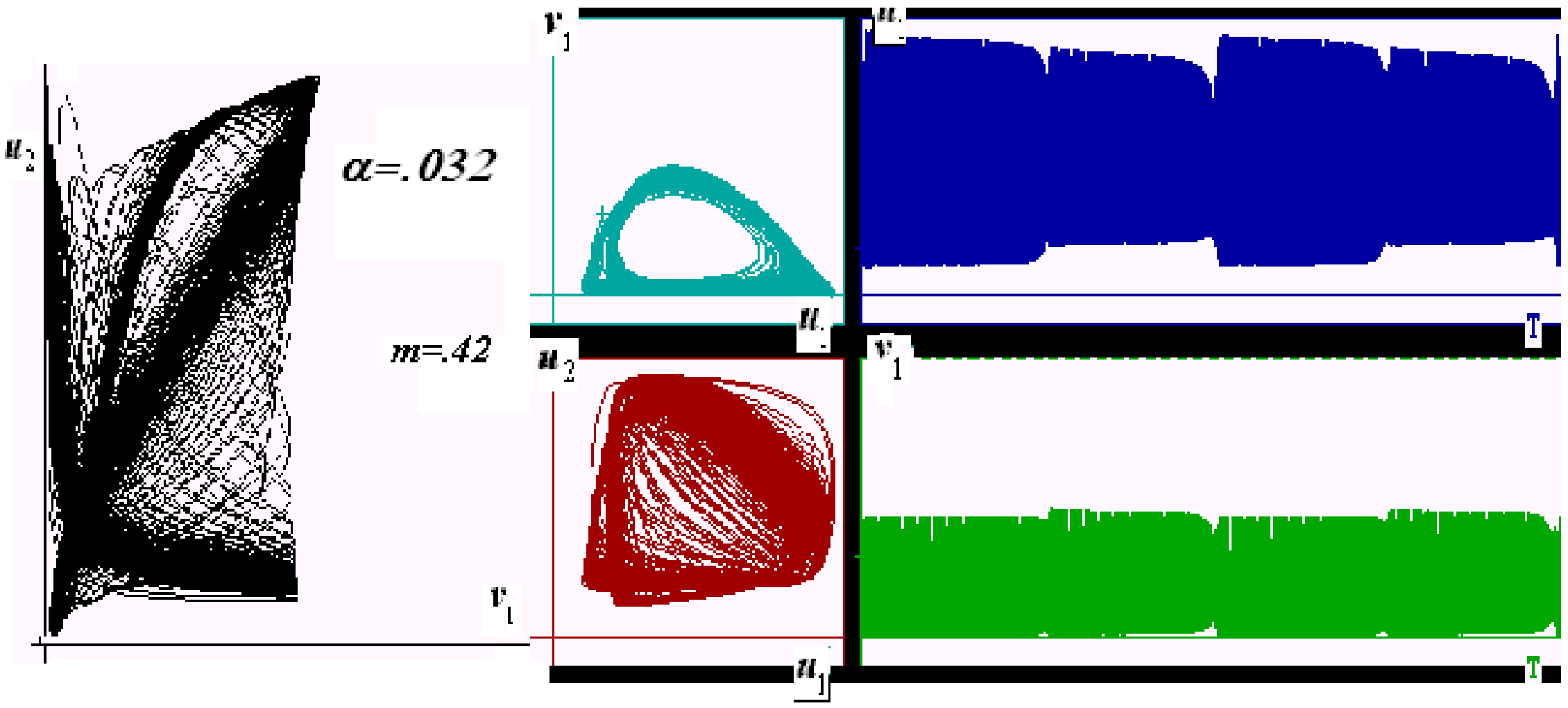

Fig.12. Irrational rotating number at $c_4$ in domain **V**

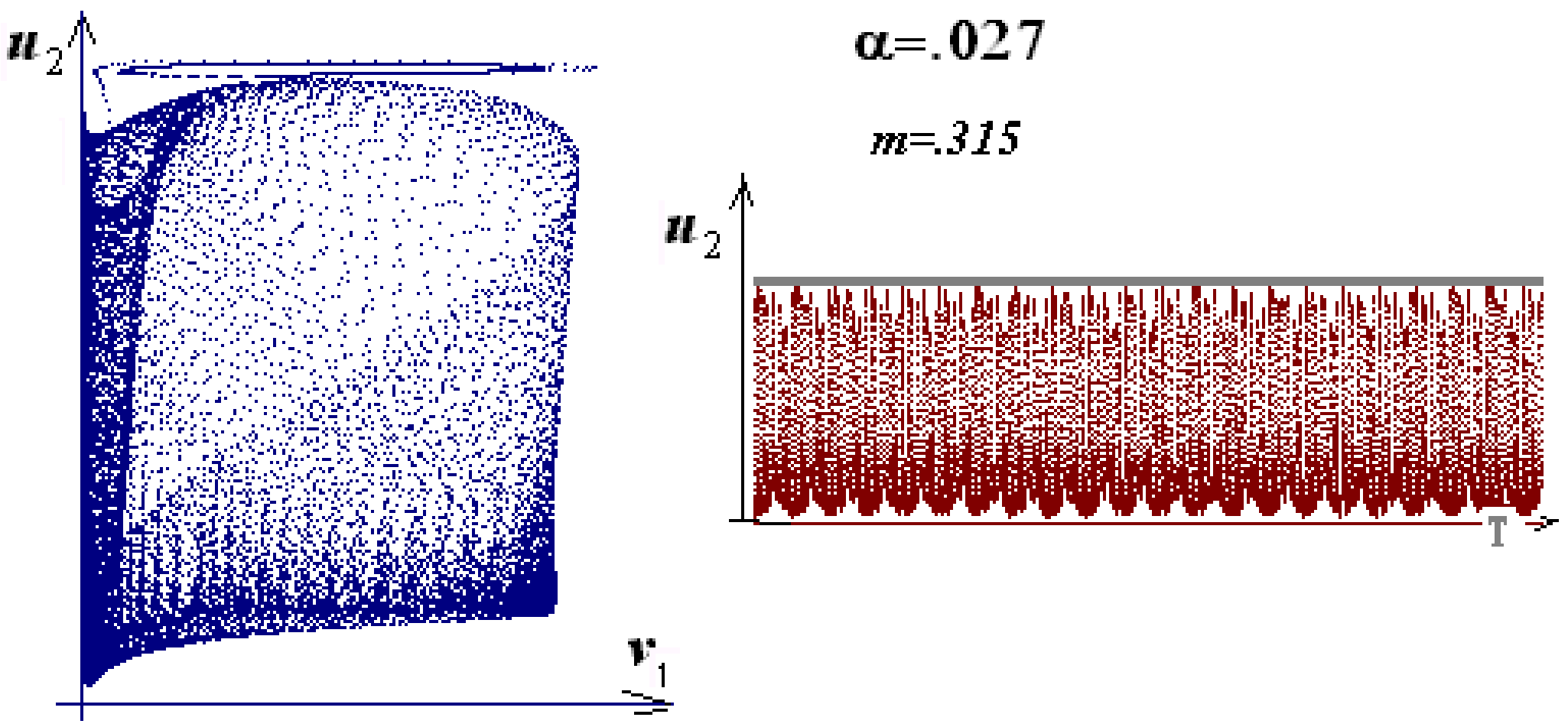

Fig.13.Stable"cycles" $c_4$ and $c_3$ coexists in domain **V**

**a** **b**

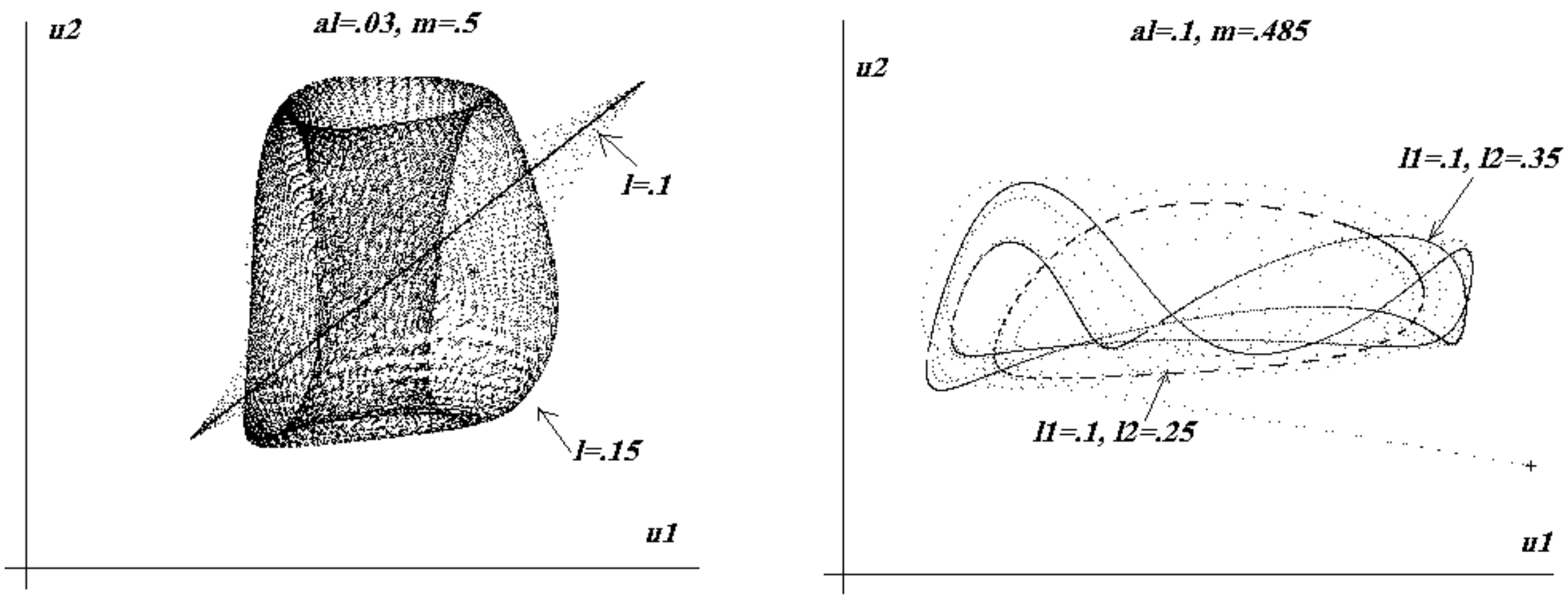

Fig.14. Examples of dynamical modes arising in the model if parameters $l_1 \neq l_2$

**Table 1**. Asymptotic coordinates and eigenvalues associated with the two-dimensional trivial equilibria $C$ of model (1s)

| Equilibrium | $\lambda_1$ | $\lambda_2$ | $\lambda_3$ | $\lambda_4$ |
|---|---|---|---|---|
| $C_{0l}(\alpha,0,l+\frac{\alpha}{1-l},0)$ | $-\gamma(m-\alpha)$ | $-\gamma(m-\frac{\alpha}{1-l})$ | $-l+(1+2l)\alpha$ | $l(1-l)+\frac{(1-3l)\alpha}{1-l}$ |
| $C_{01}(\frac{\alpha}{l},0,1-\frac{\alpha}{1-l},0)$ | $-\gamma(m-\frac{\alpha}{l})$ | $\gamma(1-m+\frac{\alpha}{1-l})$ | $-l+\frac{(2+l)\alpha}{l}$ | $-(1-l)+\frac{(3-l)\alpha}{1-l}$ |
| $C_{l0}(l+\frac{\alpha}{1-l},0,\alpha,0)$ | $-\gamma(m-l-\frac{\alpha}{1-l})$ | $-\gamma(m-\alpha)$ | $l(1-l)+\frac{(1-3l)\alpha}{1-l}$ | $-l+(1+2l)\alpha$ |
| $C_{l1}(l-\frac{\alpha}{l},0,1-\alpha,0)$ | $-\gamma(m-1+\frac{\alpha}{1-l})$ | $\gamma(1-m+\alpha)$ | $l(1-l)-\frac{(2-3l)\alpha}{l}$ | $-(1-l)+(3-2l)\alpha$ |
| $C_{10}(1-\frac{\alpha}{1-l},0,\frac{\alpha}{l},0)$ | $\gamma(1-m+\frac{\alpha}{1-l})$ | $-\gamma(m-\frac{\alpha}{l})$ | $-l+\frac{(2+l)\alpha}{l}$ | $-(1-l)+\frac{(3-l)\alpha}{1-l}$ |
| $C_{1l}(1-\alpha,0,l-\frac{\alpha}{l},0)$ | $\gamma(1-m+\alpha)$ | $-\gamma(m-l+\frac{\alpha}{l})$ | $1-l-\frac{(2-3l)\alpha}{l}$ | $-(1-l)+(3-2l)\alpha$ |

**Table 2**. Asymptotic coordinates up to $O(\alpha^2)$ and eigenvalues up to $O(\alpha)$ associated to the three-dimensional trivial equilibria $B^1(m,y,x,0)$, $B^2(x,0,m,y)$, of Model (3) for parameter values $27\alpha m - 9(\alpha + l)(1+l) + 2(1+l)^3 + 4(3\alpha - 1 + l - l^2)^3 < 0$

| | $B_0^1(m,y,x,0)$ | $B_l^1(m,y,x,0)$ | $B_1^1(m,y,x,0)$ |
|---|---|---|---|
| | $x = m\frac{\alpha}{l} + m(m-l+ml)\frac{\alpha^2}{l^3}$, $y = (m-l)(1-m) - \alpha(1-\frac{\alpha}{l})$ | $x = l - (m-l)\frac{\alpha}{l(1-l)}$, $y = (m-l)(1-m) - \frac{(m-l)}{m}\alpha$ | $x = 1 - (1-m)\frac{\alpha}{l(1-l)} - (1-m)(1-2m+lm)\frac{\alpha^2}{(1-l)^3}$, $y = (m-l)(1-m) - \frac{(1-m)\alpha}{m}$ |
| $\lambda_1$ | $-\gamma m(1-\frac{\alpha}{l})$ | $-\gamma(m-l)(1-\frac{\alpha}{l(1-l)})$ | $\gamma\frac{1-m}{1-l}(1-\frac{(1+lm-2m)\alpha}{(1-l)^2})$ |
| $\lambda_2$ | $-l + (2ml + 2m - l)\frac{\alpha}{l}$ | $l(1-l) - \frac{2(1-2l)(m-l)\alpha}{l(1-l)}$ | $-(1-l) + \frac{(3-l-4m+2lm)\alpha}{(1-l)}$ |
| $\lambda_{3,4}$ | $\frac{Tr_1 \pm \sqrt{\delta_1}}{2}$, *where* | $\frac{Tr_2 \pm \sqrt{\delta_2}}{2}$, *where* | $\frac{Tr_3 \pm \sqrt{\delta_3}}{2}$, *where* |
| | $Tr_1 = (1+l-2m)m - (1-m)\alpha$ | $Tr_2 = (1+l-2m)m - (1+l-m)\alpha$ | $Tr_3 = (1+l-2m)m - (2-m)\alpha$ |
| | $\delta_1 = Tr_1^2 - 4\gamma m(m-l)(1-m) + 4\gamma m\alpha$ | $\delta_2 = Tr_2^2 - 4\gamma m(m-l)(1-m) + 4\gamma(m-l)\alpha$ | $\delta_3 = Tr_3^2 - 4\gamma m(m-l)(1-m) - 4\gamma(1-m)\alpha$ |